\documentclass[11pt, twoside]{amsart}

\usepackage[T1]{fontenc}
\usepackage{amssymb,amsmath,amstext}

\newtheorem{theor}{Theorem}[section]
\newtheorem{lem}[theor]{Lemma}
\newtheorem{defin}[theor]{Definition}

\newtheorem{cor}[theor]{Corollary}
\newtheorem{rem}[theor]{Remark}

\newtheorem{fact}[theor]{Fact}

\newtheorem{assump}[theor]{Assumption}
\newtheorem{observation}[theor]{Observation}

\newtheorem{terminologyandnotation}[theor]{Terminology and notation}

\newcommand{\mr}{\mathrm}
\newcommand{\mc}{\mathcal}
\newcommand{\mb}{\mathbf}

\newcommand{\es}{\emptyset}

\newcommand{\nts}{\negthickspace}
\newcommand{\uhrc}{\nts \upharpoonright \nts}

\newcommand{\dist}{\mathrm{dist}}

\newcommand{\mcG}{\mathcal{G}}
\newcommand{\mcH}{\mathcal{H}}

\newcommand{\mcK}{\mathcal{K}}

\newcommand{\mcN}{\mathcal{N}}

\newcommand{\mbA}{\mathbf{A}}

\newcommand{\mbC}{\mathbf{C}}
\newcommand{\mbD}{\mathbf{D}}
\newcommand{\mbE}{\mathbf{E}}
\newcommand{\mbF}{\mathbf{F}}

\newcommand{\mbP}{\mathbf{P}}
\newcommand{\mbQ}{\mathbf{Q}}

\newcommand{\mbK}{\mathbf{K}}

\newcommand{\mbX}{\mathbf{X}}

\newcommand{\mbbN}{\mathbb{N}}
\newcommand{\mbbR}{\mathbb{R}}

\newcommand{\mbbQ}{\mathbb{Q}}

\title[Almost $l$-partite graphs]
{A limit law of almost $l$-partite graphs}

\author{Vera Koponen}

\address{Vera Koponen, Department of Mathematics, Uppsala University, Box 480,
75106 Uppsala, Sweden.}

\email{vera@math.uu.se}

\begin{document}

\begin{abstract} 
For integers $l \geq 1$, $d \geq 0$ we study (undirected) graphs with vertices 
$1, \ldots, n$ such that the vertices
can be partitioned into $l$ parts such that every vertex has at most
$d$ neighbours in its own part.
The set of all such graphs is denoted $\mbP_n(l,d)$.
We prove a labelled first-order limit law, i.e., for every first-order sentence
$\varphi$, the proportion of graphs in $\mbP_n(l,d)$ that satisfy $\varphi$ converges
as $n \to \infty$. 
By combining this result with a result of 
Hundack, Prömel and Steger \cite{HPS} we also prove that if
$1 \leq s_1 \leq \ldots \leq s_l$ are integers, then $\mb{Forb}(\mcK_{1, s_1, \ldots, s_l})$
has a labelled first-order limit law, where 
$\mb{Forb}(\mcK_{1, s_1, \ldots, s_l})$ denotes the set of all graphs with vertices $1, \ldots, n$,
for some $n$, in which there is no subgraph isomorphic to 
the complete $(l+1)$-partite graph with parts of sizes $1, s_1, \ldots, s_l$.
In the course of doing this we also prove that there exists a first-order formula
$\xi$, depending only on $l$ and $d$,
such that the proportion of $\mcG \in \mbP_n(l,d)$ with the following property
approaches 1 as $n \to \infty$: there is a unique partition of $\{1, \ldots, n\}$ into
$l$ parts such that every vertex has at most $d$ neighbours in its own part, and
this partition, viewed as an equivalence relation, is defined by $\xi$.
\medskip

\end{abstract}

\maketitle

\section{Introduction}\label{introduction}

\noindent
Over the last four decades a large number of logical limit laws and
zero-one laws, as well as some non-convergence results, have been proved,
for various collections of finite structures, various probability measures
and various logics. One of the main directions of research has considered random graphs with
vertex set $[n] = \{1, \ldots, n\}$ such that, for some $0 < \alpha < 1$, an edge appears
between two vertices with probability $n^{-\alpha}$, independently of other edges.
See \cite{SS, Sp} for this line of research.
This article deals with the following context. For a first-order language $L$
and every positive integer $n$,
let $\mbK_n$ be a set of $L$-structures with universe $[n]$,
so we are dealing with `labelled' structures.
Give all members of $\mbK_n$ the same probability $1/|\mbK_n|$, so the probability
that a random member of $\mbK_n$ belongs to $\mbX \subseteq \mbK_n$ equals the proportion $|\mbX|/|\mbK_n|$.
We say that $\mbK = \bigcup_{n \in \mbbN^+} \mbK_n$ has a {\em limit law} if
for every $L$-sentence $\varphi$, the proportion of $\mcG \in \mbK_n$ 
in which $\varphi$ is true converges as $n$ tends to infinity. 
If the limit is always~0 or~1 then we say that $\mbK$ has a {\em zero-one law}.
Such a result was first proved by Glebskii, Kogan, Liogonkii,
and Talanov \cite{Gleb} and independently by
Fagin \cite{Fag} in the case when $\mbK_n$ contains all $L$-structures with
universe $[n]$ and $L$ has finite relational vocabulary and every relation symbol has
arity at least 2.
Suppose that we keep the assumptions on the language, but restrict membership in $\mbK_n$
to $L$-structures with universe $[n]$ which satisfy some constraints.
What can we say about limit laws in this case?
In general, dividing lines for when a limit law holds, or not, are not known.
But a number of results have been obtained for various $\mbK$.
Compton \cite{Com88} has proved that if $\mbK_n$ is the set of partial orders,
then $\mbK$ satisfies a zero-one law. Compton \cite{Com87} and others have also developed a theory of
limit laws (with emphasis on `unlabelled' structures)
when $\mbK$ is, up to isomorphism, closed under forming disjoint unions
and extracting connected components and the growth of $|\mbK_n|$ is slow as $n$ grows.
A book by Burris \cite{Bur} treats this theory, based on number theory.
Kolaitis, Prömel and Rothschild \cite{KPR} have proved a zero-one law in the
case when $\mbK_n$ is the set of $(l+1)$-clique free graphs ($l \geq 2$).
In the process of doing this they proved that if $\mbK_n$ is the set of
$l$-partite (or $l$-colourable) graphs, then $\mbK$ satisfies a zero-one law.
This result was generalised by the author who proved that whenever the 
vocabulary of $L$ is finite, relational and all relation symbols have arity at least 2,
then, with $\mbK_n$ being the set of $l$-colourable $L$-structures, $\mbK$ has a
zero-one law \cite{Kop09}.
Lynch \cite{Lyn} has proved a limit law when (for every $n$) $\mbK_n$ consists of all graphs with
a degree sequence that satisfies certain conditions; in particular his result implies that
$\mbK$ has a limit law when it is the set of $d$-regular graphs ($d$ fixed)
with vertex set $[n]$ for some $n$.
More results about limit laws when $\mbK_n$ is the set of $d(n)$-regular graphs and $d(n)$ is a growing function
appear in work of Haber and Krivelevich \cite{HK}.
In the case when $\mbK_n$ is the set of graphs
with vertex set $[n]$ in which every vertex has degree at most $d$,
a limit law also holds \cite{Kop12}. That will be used in this paper.

The author has two viewpoints on the present work.
One is that it adds more examples of collections of structures for which a limit law holds.
In particular, we get more examples of graphs $\mcH$ for which the set of $\mcH$-free graphs
satisfy a limit law (but in general not a zero-one law). The only previously known
example appears to be when $\mcH$ is an $(l+1)$-clique for $l \geq 2$ \cite{KPR}.
The addition of more concrete examples may be of help in attempts to understand
dividing lines between $\mbK$ with a limit law and $\mbK$ without it.

Another viewpoint is that the work presented here seeks to develop
methods for understanding limit laws in the case when members of
$\mbK$ can be decomposed into simpler substructures (in some sense) and
where the interaction between these substructures is known (at least in a
probabilistic sense). In particular, we will use knowledge from \cite{Kop12} about the 
typical structure of graphs with maximum degree $d$ when studying
$\mbK_n = \mbP_n(l,d)$, the set of graphs with vertex set $[n]$ such that $[n]$
can be partitioned into $l$ parts such that every vertex has at most 
$d$ neighbours in its own part.
This approach to understanding asymptotic properties is inspired by infinite model
theory, where one often tries to understand structures in terms of
simpler building blocks (strongly minimal sets, rank one sets, etc.)
and how these blocks are ``glued'' together.

When proving a limit law for $\mcH$-free graphs where $\mcH$
is as in Theorem~\ref{limit law for forbidden l+1-partite graphs}, below,
we use a structure result for almost all $\mcH$-free graphs  by
Hundack, Prömel and Steger \cite{HPS}
(when $\mcH$ has a colour critical vertex, defined below).
More structural results for other choices of $\mcH$ and almost all $\mcH$-free
graphs have been proved by Balogh, Bollob\'{a}s and Simonovits \cite{BBS09, BBS11}.
These may be useful in further studies of limit laws.

We now describe the main results this paper.
By `graph' we always mean `undirected graph'.
For integers $l \geq 1$ and $d \geq 0$ let 
$\mbP_n(l,d)$ be the set of graphs with vertex set $[n] = \{1, \ldots, n\}$
such that $[n]$ can be partitioned into $l$ parts such that every vertex has
at most $d$ neighbours in its own part.
Let $\mbP(l,d) = \bigcup_{n \in \mbbN^+} \mbP_n(l,d)$.
Note that $\mbP_n(l,0)$ is the set of $l$-partite, or $l$-colourable, graphs
with vertex set $[n]$, and that $\mbP_n(1,d)$ is the set of graphs with vertex set
$[n]$ in which every vertex has degree at most $d$.
For integers $1 \leq s_1 \leq s_2 \leq \ldots \leq s_l$ let $\mcK_{1, s_1, s_2, \ldots, s_l}$ denote
the complete $(l+1)$-partite graph with parts of sizes $1, s_1, s_2, \ldots, s_l$.
So if $s_1 = \ldots = s_l = 1$ then $\mcK_{1, s_1, s_2, \ldots, s_l}$ 
is an $(l+1)$-clique, i.e. a complete graph on $l+1$ vertices.
For a graph $\mcH$ let $\mb{Forb}_n(\mcH)$ be the set of graphs with vertex
set $[n]$ which contain no subgraph that is isomorphic to $\mcH$,
and let $\mb{Forb}(\mcH) = \bigcup_{n \in \mbbN^+} \mb{Forb}_n(\mcH)$.
Note that $\mbP_n(2,0) \subseteq \mb{Forb}_n(\mcK_{1,1,1})$.
In an article from 1976 \cite{EKR}, Erdös, Kleitman and Rothschild proved
that the proportion of $\mcG \in \mb{Forb}_n(\mcK_{1,1,1})$ which
are bipartite, i.e., belong to $\mbP_n(2,0)$, approaches 1 as $n \to \infty$.
Later, Kolaitis, Prömel and Rothschild \cite{KPR} generalised this by
proving that, for every $l \geq 2$, if $s_1 = s_2 = \ldots = s_l = 1$, then
$|\mbP_n(l,0)| \big/ |\mb{Forb}_n(\mcK_{1, s_1, \ldots, s_l})| \to 1$ as $n \to \infty$ and
$\mbP(l,0)$ satisfies a zero-one law; hence also $\mb{Forb}_n(\mcK_{1, s_1, \ldots, s_l})$
satisfies a zero-one law if $s_1 = s_2 = \ldots = s_l = 1$.

We say that a vertex $v$ of a graph $\mcH$ is {\em colour-critical} if 
one can obtain a graph with smaller chromatic number than $\mcH$ by removing some
edges of $\mcH$ which contain $v$, and only such edges.
The {\em criticality} of a colour-critical vertex $v$ is the minimal 
number of edges which contain $v$ that must be removed to produce a graph with smaller 
chromatic number.
Prömel and Steger \cite{PS} and then Hundack, Prömel and Steger \cite{HPS} have generalised
the result of Kolaitis, Prömel and Rothschild that almost all $(l+1)$-clique-free graphs
are $l$-partite to the following:

\begin{theor}\label{HPS-results} \cite{HPS}
Suppose that $\mcH$ is a graph with chromatic number $l+1$
and with a colour critical vertex $v$ 
with criticality $d$ and suppose that no other colour-critical vertex has smaller criticality than $v$.
Then 
$$\frac{|\mb{Forb}_n(\mcH) \cap \mbP_n(l,d-1)|}{|\mb{Forb}_n(\mcH)|} \ \to 1 \ \text{ as } n \to \infty.$$
\end{theor}

\noindent
The main result of this article is the following, where the `language of graphs' refers to
the first-order language built up from a vocabulary (also called signature) which consists
only of a binary relation symbol, besides the identity symbol:

\begin{theor}\label{limit theorem}
Suppose that $l \geq 1$ and $d \geq 0$ are integers.
For every first-order sentence $\varphi$ in the language 
of graphs, the proportion of $\mcG \in \mbP_n(l,d)$ in which $\varphi$ is true
converges as $n \to \infty$. If $d = 0$ or $d = 1$ then the this proportion
always converges to either 0 or 1; if $d > 1$ then it may converge to some $0 < c < 1$.
\end{theor}

\noindent
In the case $d = 0$ Theorem~\ref{limit theorem} states the same thing as
one of the main results of \cite{KPR} (described above). In the case $l = 1$
Theorem~\ref{limit theorem} states the same thing as the main result of \cite{Kop12}.
Therefore we focus on the case when $d \geq 1$ and $l \geq 2$.
Theorems~\ref{HPS-results} and~\ref{limit theorem} will be used to prove
the following result, in the last section.

\begin{theor}\label{limit law for forbidden l+1-partite graphs}
Suppose that $l \geq 2$, $1 \leq s_1 \leq s_2 \leq \ldots \leq s_l$ are integers.\\
(i) For every sentence $\varphi$ in the language of graphs, the 
proportion of $\mcG \in \mb{Forb}_n(\mcK_{1,s_1, \ldots, s_l})$ in which $\varphi$ is true
converges as $n \to \infty$. \\
(ii) If $s_1 \leq 2$ then this proportion converges to 0 or 1 for every sentence $\varphi$.\\
(iii) If $s_1 > 2$ then there are infinitely many mutually contradictory sentences 
$\varphi_i$, $i \in \mbbN$, in the language of graphs such that the proportion of
$\mcG \in \mb{Forb}_n(\mcK_{1,s_1, \ldots, s_l})$ in which $\varphi_i$ is true approaches
some $\alpha_i$ such that $0 < \alpha_i < 1$.
\end{theor}

\noindent
This article is organised as follows.
Section~\ref{decompositions} considers the possibly different
ways in which the vertex set of $\mcG \in \mbP_n(l,d)$ can be partitioned
into $l$ parts such that every vertex has at most $d$ neighbours
in its own part. 
We show that there is $\mu > 0$, depending only on $l$, 
such that the proportion of $\mcG \in \mbP_n(l,d)$ with the
following property approaches 1 as $n\to\infty$:
for every partition $V_1, \ldots, V_l$ of the vertex set such that every
vertex has at most $d$ neighbours in its own part, $|V_i| \geq \mu n$ for all $i \in [l]$.
In Section~\ref{extension properties} we consider the following sort of question,
the probability of an ``extension property'',
where $\mcH_1$ is assumed to be an induced subgraph of $\mcH_2$:
Given $\mcG \in \mbP_n(l,d)$, what is the probability that
every induced subgraph of $\mcG$ that is isomorphic to $\mcH_1$ is contained
in an induced subgraph of $\mcG$ which is isomorphic to $\mcH_2$?
In Section~\ref{unique decomposition}
we use the results from sections~\ref{decompositions} and~\ref{extension properties}
to prove that the proportion of $\mcG \in \mbP_n(l,d)$ with the following 
property approaches 1 as $n \to \infty$:
there is exactly one way in which the vertex set can be partitioned
into $l$ (non-empty) parts such that every vertex has at most $d$ neighbours
in its own part.
In Section~\ref{a limit law}
we use the results from Sections~\ref{extension properties} and~\ref{unique decomposition},
the main results from \cite{Kop12} and an Ehrenfeucht-Fra\"{i}ss\'{e} game argument to prove
Theorem~\ref{limit theorem}.
In Section~\ref{forbidden subgraphs} we consider ``forbidden subgraphs''
of the type $\mcK_{1, s_1, \ldots, s_l}$ and prove 
Theorem~\ref{limit law for forbidden l+1-partite graphs}
with the help of Theorem~\ref{limit theorem}.

\begin{terminologyandnotation}{\rm
See for example \cite{EF} for an introduction to first-order logic and first-order structures
and \cite{Die} for basics about graph theory.
By {\em graph} we mean undirected graph without loops.
By the {\em first-order language of graphs} we mean the 
set of first-order formulas over a vocabulary (also called signature) with 
the identity symbol `=' and a binary relation symbol `$E$'
(for the edge relation). When speaking of a formula or sentence we
will always mean a formula, or sentence, in the language of graphs.
We view graphs as first-order structures $\mcG = (V, E^{\mcG})$ for the language of graphs.
Since we only consider undirected graphs without loops, the interpretation
of $E$, $E^{\mcG}$, will always be symmetric and irreflexive, so we may, if convenient, 
view $E^{\mcG}$ as a set of 2-subsets of $V$.
Let $\mcG = (V, E^{\mcG})$ be a graph.
If $\varphi(x_1, \ldots, x_m)$ is a formula with free variables $x_1, \ldots, x_m$
and $v_1, \ldots, v_m \in V$, then the notation `$\mcG \models \varphi(v_1, \ldots, v_m)$' 
means that $v_1, \ldots, v_m$ satisfies the statement $\varphi(x_1, \ldots, x_m)$ in $\mcG$,
and for a sentence $\psi$, $\mcG \models \psi$ means that $\psi$ is satisfied by $\mcG$
(or in other words, that $\mcG$ has the property expressed by $\psi$).
For $v, w \in V$, the notation $v \sim_{\mcG} w$ means that $v$ and $w$ are adjacent in $\mcG$;
so $v \sim_{\mcG} w$ expresses the same thing as $\mcG \models E(v,w)$.
We say that $\mcH = (W, E^{\mcH})$ is a {\em subgraph} of $\mcG = (V, E^{\mcG})$
if $W \subseteq V$ and $E^{\mcH} \subseteq E^{\mcG}$.
If, in addition, for all $a \in W$, $a \sim_{\mcH} b$ if and only if $a \sim_{\mcG} b$,
then we call $\mcH$ an {\em induced subgraph} of $\mcG$.
Hence, $\mcH$ is an induced subgraph of $\mcG$ if and only if $\mcH$ is a substructure of $\mcG$
in the sense of model theory.
For $X \subseteq V$, $\mcG[X]$ denotes the induced subgraph of $\mcG$ with vertex set $X$.
In model theoretic terms, $\mcG[X]$ is the substructure of $\mcG$ with universe $X$.
The distance between two vertices $v, w \in V$ in $\mcG$ is denoted $\mathrm{dist}_{\mcG}(v,w)$,
and for sets of vertices $A$ and $B$, 
$\dist_{\mcG}(A,B) = \min \{\dist_{\mcG}(v,w) : v \in A, \ w \in B\}$.
If $\mcG = (V, E^{\mcG})$ and $\mcH = (W, E^{\mcH})$ are graphs and 
$f : V \to W$ is injective and has the property that, for all $a, b \in V$,
$a \sim_{\mcG} b$ if and only if $f(a) \sim_{\mcH} f(b)$, then we call
$f$ a {\em strong embedding} of $\mcG$ into $\mcH$.
We say that functions $f,g : \mbbN \to \mbbR$ are {\em asymptotic}, written $f \sim g$,
if $f(n)/g(n) \to 1$ as $n \to \infty$.
}\end{terminologyandnotation}

\section{Decompositions}\label{decompositions}

\noindent
Let $l \geq 1$ and $d \geq 0$ be integers.
Let $\mbP_n(l,d)$ be the set of graphs with vertex set $[n] = \{1, \ldots, n\}$
such that $[n]$ can be partitioned into $l$ parts in such a way
that every vertex has at most $d$ neighbours in its own part.
In general, for $\mcG \in \mbP_n(l,d)$ there may be more than one 
partition of the vertex set into $l$ parts such that every vertex
has at most $d$ neighbours in its own part.
In this section we show that there is $\mu > 0$ depending only on $l$
such that for almost all $\mcG \in \mbP_n(l,d)$ (for large enough $n$) every
such partition $V_1, \ldots, V_l$ has the property that $|V_i| \geq \mu n$
for all $i = 1, \ldots, l$.

\begin{defin}\label{definition of decompositions}{\rm
Let $\mcG = (V,E^{\mcG}) \in \mbP_n(l,d)$, so $V = [n]$.
By the definition of $\mbP_n(l,d)$,
there exists a partition $V_1, \ldots, V_l$ of $V$, which we denote by $\pi$,
such that the following holds:
\begin{itemize}
\item[(i)] $E^{\mcG} = E_1 \cup E_2$,
\item[(ii)] The graph $(V,E_1)$ is $l$-colourable and 
the partition $\pi$ defines an $l$-colouring of it.
\item[(iii)] $E_2 = E'_1 \cup \ldots \cup E'_l$ and for every $i \in [l]$, 
$E'_i \subseteq V_i^{(2)}$ and every vertex of the graph $(V_i,E'_i)$
has degree $\leq d$.
\end{itemize}
A pair $(E_1, E_2)$ such that (i)--(iii) hold
is called a {\em decomposition of $\mcG$ based on $\pi$}, or
{\em a $\pi$-based decomposition of $\mcG$}. 
A pair $(E_1, E_2)$ is called a {\em decomposition of $\mcG$} if, for some partition $\pi$ of $V$
into $l$ parts, 
it is a $\pi$-based decomposition of $\mcG$.
}\end{defin}

\noindent
Partitions of $V = [n]$ into $l$ parts will be denoted by $\pi$, sometimes with an index.
Note that if $\mcG = (V, E^{\mcG}) \in \mbP_n(l,0)$ and $(E_1, E_2)$ is a decomposition of
$\mcG$ which is based on a partition $\pi$ of $V$ into $l$ parts $V_1, \ldots, V_l$, then
$E_2 = \es$, $E_1 = E^{\mcG}$ and
$\pi$ induces an $l$-colouring of $\mcG$, 
in the sense that all elements in $V_i$ can be assigned the colour $i$, for all $i \in [l]$.
It is straightforward to verify the following:

\begin{observation}\label{observation about decompositions}{\rm
Let $\mcG = (V,E^{\mcG}) \in \mbP_n(l,d)$.\\
(a) If $(E_1, E_2)$ is a decomposition of $\mcG$, 
then $E_1$ and $E_2$ are disjoint.\\
(b) By the definition of $\mbP_n(l,d)$, in the beginning
of the section, $\mcG$ has a
decomposition based on some partition of $V = [n]$ into $l$ parts.\\
(c) For every partition $\pi$ of $V$ into $l$ parts, there is at most
one decomposition of $\mcG$ which is based on $\pi$.\\
(d) In general, it is possible that there are different partitions of $V$ into $l$ parts, 
say $\pi_1$ and $\pi_2$,
a decomposition of $\mcG$ based on $\pi_1$ and another decomposition of $\mcG$
based on $\pi_2$.
}\end{observation}

\noindent
Part (d) of the observation might look discouraging because, in general,
there is not a unique way to present a graph $\mcG \in \mbP(l,d)$ by its decomposition.
However, the next couple of lemmas together with the results in
Sections~\ref{extension properties}, \ref{unique decomposition}
show that, as $n \to \infty$,
the proportion of graphs $\mcG \in \mbP_n(l,d)$ which have a unique decomposition
approaches~1.

\begin{defin}\label{definition of richness}{\rm
(i) Let $\alpha \in \mbbR$. An $l$-colouring $f : [n] \to [l]$ of a graph $\mcG \in \mbP_n(l,0)$ is called
{\em $\alpha$-rich} if $|f^{-1}(i)| \geq \alpha$ for every $i \in [l]$; that is, for every colour $i$,
at least $\alpha$ vertices are assigned the colour $i$ by $f$.\\
(ii) Similarly as in (i), a partition of $[n]$ into $l$ parts is called {\em $\alpha$-rich}
if each one of the $l$ parts contains at least $\alpha$ elements.\\
(iii) Let $\pi$ denote any partition of $[n]$ into $l$ parts.
By $\mbP_{n,\pi}(l,d)$ we denote the set of all $\mcG \in \mbP_n(l,d)$
which have a $\pi$-based decomposition.
}\end{defin}

\noindent
Theorem~10.5 in \cite{Kop09} has the following as an immediate consequence:

\begin{fact}\label{fast convergence of mu-n-rich colourings}\cite{Kop09}
For every $l \geq 2$ and every sufficiently small $\mu > 0$, there is $\lambda > 0$ such that for 
all sufficiently large $n$,
\begin{equation}\label{formula for fast convergence}
\frac{\big|\{\mcG \in \mbP_n(l,0) : \text{ $\mcG$ has an $l$-colouring which is {\em not}
$\mu n$-rich}\}\big|}{|\mbP_n(l,0)|} \ \leq \ 2^{- \lambda n^2 \pm O(n)}.
\end{equation}
\end{fact}

\noindent
An analysis of the proof of Theorem~10.5~(ii) in~\cite{Kop09} shows that if
$0 < \mu < \frac{1}{2l(l-1)}$, then there is $\lambda > 0$ such that the conclusion of
Fact~\ref{fast convergence of mu-n-rich colourings} holds.
By applying Fact~\ref{fast convergence of mu-n-rich colourings} we get the following:

\begin{cor}\label{almost all graphs have a partition with large parts}
Let $l \geq 1$ and $d \geq 0$ be integers.
If $\mu > 0$ is sufficiently small and $\widehat{\mbP}_n(l,d)$
denotes the set of all $\mcG \in \mbP_n(l,d)$ which have
a decomposition which is based on a partition 
that is {\em not} $\mu n$-rich, then there is $\lambda > 0$ (depending on $\mu$)
such that for all sufficiently large $n$,
\begin{equation*}
\frac{\big|\widehat{\mbP}_n(l,d)\big|}{\big|\mbP_n(l,d)\big|} \ \leq \
2^{-\lambda n^2 + O(n \log n)} \ \to \ 0 \ \ \text{ as } n \to \infty.
\end{equation*}
\end{cor}

\noindent
{\bf Proof.}
If $l = 1$ then for every $0 < \mu \leq 1$ we have $\widehat{\mbP}_n(l,d) = \es$ for all $n$,
so the conclusion of the lemma follows trivially.
Now suppose that $l > 1$.
By Fact~\ref{fast convergence of mu-n-rich colourings}, we can choose 
$\mu > 0$ small enough so that there exists $\lambda > 0$ 
such that~(\ref{formula for fast convergence}) holds
for all sufficiently large $n$.
Recall that $\mbP_n(l,0)$ is the set of all $l$-colourable graphs with vertices $1, \ldots, n$.
Also, observe that
$\mbP_n(1,d)$ is the set of all graphs with vertices $1, \ldots, n$
such that every vertex has degree $\leq d$.
Note that if $V = [n]$ and 
$(E_1, E_2)$ is a decomposition of $\mcG \in \mbP_n(l,d)$,
then $(V, E_1) \in \mbP_n(l,0)$ and 
$(V,E_2) \in \mbP_n(1,d)$.
Observe that $\mcG \in \mbP_n(l,0)$ has an $l$-colouring which is {\em not}
$\mu n$-rich if and only if $\mcG$ has a decomposition based on a partition
which is {\em not} $\mu n$-rich.
Since every $\mcG \in \widehat{\mbP}_n(l,d)$ has 
a decomposition which is based on a partition of $V = [n]$ into $l$ parts
which is {\em not} $\mu n$-rich, it follows that
\begin{equation}\label{description in terms of colouring and bounded degree}
\big|\widehat{\mbP}_n(l,d)\big| \ \leq \ 
\big|\widehat{\mbP}_n(l,0)\big| \cdot \big|\mbP_n(1,d)\big|,
\end{equation}
where $\widehat{\mbP}_n(l,0)$ is the set of $\mcG \in \mbP_n(l,0)$
that have an $l$-colouring which is not $\mu n$-rich. 
Next, we estimate an upper bound of $|\mbP_n(1,d)|$.
For all sufficiently large $n$, each vertex of a graph in $\mbP_n(1,d)$
can be connected to the other vertices in at most
$\sum_{i=0}^{d-1} \binom{n}{i} \leq  dn^{d-1}$ ways. 
Therefore,
\begin{equation}\label{upper bound of graphs with bounded degree}
\big|\mbP_n(1,d)\big| \ \leq \ \Big(dn^{d-1}\Big)^n \ \leq \
2^{(d-1)n\log n \ + \ O(n)}.
\end{equation}
Hence, for all sufficiently large $n$,
\begin{align*}
\frac{\big|\widehat{\mbP}_n(l,d)\big|}{\big|\mbP_n(l,d)\big|} \ &\leq \ 
\frac{\big|\widehat{\mbP}_n(l,0)\big| \cdot \big|\mbP_n(1,d)\big|}
{\big|\mbP_n(l,d)\big|} 
\quad \text{ by (\ref{description in terms of colouring and bounded degree})}\\
&\leq \
\frac{\big|\widehat{\mbP}_n(l,0)\big| \cdot \big|\mbP_n(1,d)\big|}
{\big|\mbP_n(l,0)\big|}
\quad \text{ because } \mbP_n(l,0) \subseteq \mbP_n(l,d)\\
&\leq \ 
2^{- \lambda n^2 \ \pm \ O(n)} \ \cdot \ \big|\mbP_n(1,d)\big| 
\quad \text{ by (\ref{formula for fast convergence})}\\
&\leq \
2^{- \lambda n^2 \ + \ (d-1)n \log n \ + \ O(n)} 
\quad \text{ by (\ref{upper bound of graphs with bounded degree})}\\
&= \
2^{- \lambda n^2 \ + \ O(n \log n)}. \quad \square
\end{align*}

\section{Extension properties}\label{extension properties}

\noindent
Fix an integer $d \geq 0$.
In this section we prove some technical lemmas about extension properties
which will be used in Sections~\ref{unique decomposition}
and~\ref{a limit law}.

\begin{assump}\label{assumption about H}
{\rm {\bf (until Definition~\ref{definition of k-extension property})}
Let $\mcH = (X_1 \cup X_2 \cup Y, E^{\mcH})$ be a graph where $X_1, X_2$ and $Y$ are mutually disjoint
and $v \not\sim_{\mcH} w$ whenever $v \in X_1$ and $w \in Y$.
}\end{assump}

\begin{lem}\label{first extension lemma}
For $i = 1,2$, let $\mcK_i = (W_i, E^{\mcK_i})$ be graphs 
such that $W_1 \cup W_2 = V = [n]$, $\mcK_1$ has maximum degree at most $d$
and $W_1$ has at least $n^{1/4}$ subsets $Z_i$, $i = 1, \ldots, n^{1/4}$, such that
$\mcK_1[Z_i] \cong \mcH[Y]$. 
Let $\mbP(\mcK_1, \mcK_2)$ be the set of graphs $\mcG = (V, E^{\mcG})$ such that 
$\mcG[W_i] = \mcK_i$ for $i = 1,2$.
Then the proportion of $\mcG \in \mbP(\mcK_1, \mcK_2)$ such that every
strong embedding $h_0 : \mcH[X_1 \cup X_2] \to \mcG$ satisfying $h_0(X_i) \subseteq W_i$, for $i = 1,2$,
can be extended to a strong embedding $h : \mcH \to \mcG$ satisfying $h(Y) \subseteq W_1$, is at least
$$1 - n^{|X_1|+|X_2|} \ \alpha^{\beta n^{1/4}},$$
where $0 < \alpha, \beta < 1$ are constants that depend only on $|X_1|, |X_2|, |Y|$
and $d$. 
\end{lem}

\noindent
{\bf Proof.}
Assume that $h_0 : X_1 \cup X_2 \to V = [n]$ is injective and
$h_0(X_i) \subseteq W_i$ for $i = 1,2$.
Let $\mbP_{h_0}$ be the set of $\mcG \in \mbP(\mcK_1, \mcK_2)$ such that
$h_0$ is a strong embedding of $\mcH[X_1 \cup X_2]$ into $\mcG$.
Then:
\\

\noindent
{\em Claim.} The proportion of $\mcG \in \mbP_{h_0}$ such that $h_0$ cannot
be extended to a strong embedding $h : \mcH \to \mcG$ such that $f(Y) \subseteq W_1$
is at most $(1 - 2^{-p})^{\beta n^{1/4}}$ where $p \geq 1$ and $0 < \beta < 1$
depend only on $|X_1|, |X_2|, |Y|$ and $d$.
\\

\noindent
{\em Proof of the claim.}
By one of the assumptions of the lemma, 
for every $i = 1, \ldots, n^{1/4}$, $h_0$ can be extended to a
function $h_i : X_1 \cup X_2 \cup Y \to V$ such that $h_i \uhrc Y$ is an isomorphism
from $\mcH[Y]$ onto $\mcK_1[Z_i]$.
Since the maximum degree of $\mcK_1$ is at most $d$, there are different
$i_1, \ldots, i_{\beta n^{1/4}} \in \{1, \ldots, n^{1/4}\}$, where $0 < \beta < 1$
depends only on $|X_1| + |X_2| + |Y|$ and $d$, such that 
$\dist_{\mcK_1}(h_0(X_1 \cup X_2), Z_{i_j}) > 1$ for all $j$ and
$\dist_{\mcK_1}(Z_{i_j}, Z_{i_{j'}}) > 1$ whenever $j \neq j'$.
By the assumptions about $\mcH$ (Assumption~\ref{assumption about H}), $h_{i_j} \uhrc X_1 \cup Y$ is a 
strong embedding of $\mcH[X_1 \cup Y]$ into $\mcK_1$, for every $j = 1, \ldots \beta n^{1/4}$.
From the definitions of $\mbP(\mcK_1, \mcK_2)$ and $\mbP_{h_0}$ it follows that 
for every $\mcG \in \mbP_{h_0}$ and every $j = 1, \ldots, \beta n^{1/4}$,
$h_{i_j} \uhrc X_1 \cup Y$ is a strong embedding of $\mcH[X_1 \cup Y]$ into $\mcG$.

Recall that there are no restrictions on the existence (or nonexistence)
of edges going between $W_1$ and $W_2$;
in other words we can see each such edge as existing with probability $1/2$ independently of
the other edges of the graph. 
Therefore, if all $\mcG \in \mbP_{h_0}$ have the same probability and such $\mcG$ is
chosen at random, then the probability that $h_{i_j}$ is a strong
embedding of $\mcH$ into $\mcG$ is $2^{-p}$ where $p = |X_2| \cdot |Y|$; and this
holds independently of whether $h_{i_{j'}}$ is a strong embedding of $\mcH$ into $\mcG$ for $j' \neq j$.
It follows that the proportion of $\mcG \in \mbP_{h_0}$  such that, for every $j = 1, \ldots, \beta n^{1/4}$,
$h_{i_j}$ is not a strong embedding of $\mcH$ into $\mcG$ is $\big(1 - 2^{-p}\big)^{\beta n^{1/4}}$.
\hfill $\square$
\\

\noindent
There are not more than $n^{|X_1| + |X_2|}$ injective functions 
from $X_1 \cup X_2$ into $V$. By the claim, for every such function, say $h_0$, 
the proportion of $\mcG \in \mbP_{h_0}$ such that $h_0$ cannot be extended to a strong 
embedding $h$ of $\mcH$ into $\mcG$ that satisfies $h(Y) \subseteq W_1$ is at most
$\big(1 - 2^{-p}\big)^{\beta n^{1/4}}$.
Therefore the proportion of $\mcG \in \mbP(\mcK_1, \mcK_2)$ such that there is a
strong embedding $h_0 : \mcH[X_1 \cup X_2] \to \mcG$ with $h_0(X_i) \subseteq W_i$, for $i = 1,2$,
which cannot be extended to a strong embedding $h : \mcH \to \mcG$ with $h(Y) \subseteq W_1$ is
at most $n^{|X_1| + |X_2|} \ \big(1 - 2^{-p}\big)^{\beta n^{1/4}}$.
\hfill $\square$

\begin{assump}\label{assumption on mu}{\rm
For the rest of this section we fix $\mu > 0$ small enough so that there exists $\lambda > 0$ such 
that~(\ref{formula for fast convergence}) holds for all sufficiently large $n$
and let $\pi$ denote an arbitrary $\mu n$-rich partition of $V = [n]$ into parts $V_1, \ldots, V_l$.
}\end{assump}

\noindent
Recall Definition~\ref{definition of richness}~(iii) of $\mbP_{n,\pi}(l,d)$,
for a partition $\pi$ of $[n]$.

\begin{lem}\label{second extension lemma}
Let $p \in [l]$.
Then there are constants $0 < \alpha, \beta < 1$, depending only on 
$|X_1|, |X_2|, |Y|$ and $d$, such that the proportion of $\mcG \in \mbP_{n,\pi}(l,d)$ that satisfies the following condition
is at least $1 - n^{|X_1|+|X_2|} \ \alpha^{\beta n^{1/4}}$:
\begin{itemize}
\item[$(*)$] If $\mcG[V_p]$ has at least $n^{1/4}$ different induced subgraphs which are isomorphic to $\mcH[Y]$
and $h_0 : \mcH[X_1 \cup X_2] \to \mcG$ is a strong embedding such that $h_0(X_1) \subseteq V_p$
and $h_0(X_2) \subseteq V \setminus V_p$,
then $h_0$ can be extended to a strong embedding $h : \mcH \to \mcG$ such that $h(Y) \subseteq V_p$.
\end{itemize}
\end{lem}

\noindent
{\bf Proof.}
We will reduce the proof to an application of Lemma~\ref{first extension lemma}.
Let $W_1 = V_p$ and $W_2 = V \setminus V_p$.
Suppose that $\mcK_1$ is a graph with vertex set $W_1$, maximum degree at most $d$ and such that
$\mcK_1$ has at least $n^{1/4}$ different induced subgraphs which are isomorphic to $\mcH[Y]$.
Also suppose that $\mcK_2$ is a graph such that $\mcK_2 = \mcG[W_2]$ for some $\mcG \in \mbP_{n,\pi}(l,d)$.
Then let $\mbP(\mcK_1, \mcK_2)$ be the set of graphs $\mcG$
such that $\mcG[W_1] = \mcK_1$ and $\mcG[W_2] = \mcK_2$.
Finally let $\mbQ$ be the set of $\mcG \in \mbP_{n,\pi}(l,d)$ such that $\mcG[W_1]$ (where $W_1 = V_p$)
has less than  $n^{1/4}$ different induced subgraphs which are isomorphic to $\mcH[Y]$.
Then we have
$$\mbP_{n,\pi}(l,R) \ = \ \mbQ \ \cup \ \bigcup_{\mcK_1, \mcK_2} \mbP(\mcK_1, \mcK_2),$$
where the union ranges over all pairs $(\mcK_1, \mcK_2)$ of graphs as described above.
Moreover, if $(\mcK_1, \mcK_2) \neq (\mcK'_1, \mcK'_2)$, then
$\mbP(\mcK_1, \mcK_2) \ \cap \ \mbP(\mcK'_1, \mcK'_2) = \es$.
We also have $\mbQ \ \cap \ \mbP(\mcK_1, \mcK_2) = \es$ for all pairs $(\mcK_1, \mcK_2)$ as described.
Hence we get 
$$|\mbP_{n,\pi}(l,d)| = |\mbQ| + \sum_{\mcK_1, \mcK_2} |\mbP(\mcK_1, \mcK_2)|$$
where the sum
ranges over all pairs $(\mcK_1, \mcK_2)$ of graphs as described above.
Therefore it suffices to prove that there are constants $0 < \alpha, \beta < 1$
depending only on $|X_1|$, $|X_2|$, $|Y|$ and $d$ such that
\begin{itemize}
\item[(a)] The proportion of $\mcG \in \mbQ$ which satisfy $(*)$ is at least 
$1 - n^{|X_1|+|X_2|} \ \alpha^{\beta n^{1/4}}$.
\item[(b)] For every pair $(\mcK_1, \mcK_2)$ as described above, the proportion of
$\mcG \in \mbP(\mcK_1, \mcK_2)$ which satisfy $(*)$ is at least
$1 - n^{|X_1|+|X_2|} \ \alpha^{\beta n^{1/4}}$
\end{itemize}
But (a) is trivially true because every $\mcG \in \mbQ$ has less than $n^{1/4}$ different induced 
subgraphs which are isomorphic to $\mcH[Y]$; hence $(*)$ holds for every $\mcG \in \mbQ$.
And (b) is obtained by an application of
Lemma~\ref{first extension lemma} to every pair $(\mcK_1, \mcK_2)$ as described above.
\hfill $\square$

\begin{defin}\label{definition of k-extension property}{\rm
(i) If Assumption~\ref{assumption about H} holds and
$|X_1| + |X_2| + |Y| \leq k$, then, {\em for every} $p \in [l]$, we call the condition $(*)$ of 
Lemma~\ref{second extension lemma} a {\em $k$-extension property with respect to~$\pi$}.
Note that for every $k \in \mbbN$, there are only finitely many (non-equivalent) $k$-extension properties
with respect to $\pi$.\\
(ii) We say that a graph $\mcG \in \mbP_{n,\pi}(l,d)$ {\em has the $k$-extension property} if it satisfies
every $k$-extension property with respect to $\pi$.
}\end{defin}

\begin{cor}\label{corollary to second extension property}
For every $k \in \mbbN$, there is $\varepsilon_k : \mbbN \to \mbbR$,
depending only on $k$, such that
$\lim_{n\to\infty} \varepsilon_k(n) = 0$ and the
proportion of $\mcG \in \mbP_{n,\pi}(l,d)$ which have the $k$-extension property
with respect to $\pi$
is at least $1 - \varepsilon_k(n)$.
\end{cor}

\noindent
{\bf Proof.}
There are only finitely many, say $m$, non-equivalent $k$-extension properties.
For each of these (and large enough $n$), the proportion of $\mcG \in \mbP_{n,\pi}(l,d)$ which does
{\em not} have it is,
by Lemma~\ref{second extension lemma}, at most
$n^{2k} \ \alpha^{\beta n^{1/4}}$ where $0 < \alpha, \beta < 1$
depend only on $k$ and $d$, 
so $n^{2k} \ \alpha^{\beta n^{1/4}} \to 0$ as $n \to \infty$.
Hence $\varepsilon_k(n)$ can be taken as the sum of $m$ terms of the form
$n^{2k} \ \alpha^{\beta n^{1/4}}$. 
\hfill $\square$

\noindent
The next lemma will be used in Section~\ref{a limit law}.

\begin{lem}\label{removing some edges from a graph with large extension property}
Suppose that $\mcG \in \mbP_{n,\pi}(l,d)$ and let $\mcG'$ be the graph that results
from removing (from $\mcG$) or adding (to $\mcG$) at most $s$ edges $\{v,w\}$
such that for some $i \neq j$, $v \in V_i$ and $w \in V_j$.
If $\mcG$ has the $(k + 2s)$-extension property with respect to $\pi$, 
then $\mcG'$ has the $k$-extension property with respect to $\pi$.
\end{lem}

\noindent
{\bf Proof.} Straightforward consequence of the definition of $k$-extension property. 
\hfill $\square$

\section{Unique decomposition:\\ expressing a partition in the language of graphs}
\label{unique decomposition}

\noindent
The goal in this section is to show that the proportion of $\mcG \in \mbP_n(l,d)$ which
have a unique decomposition approaches 1 as $n \to \infty$.
For $l = 1$ this is trivially true, so {\em we assume that $l \geq 2$ in this section.}
As in Assumption~\ref{assumption on mu}, 
we fix a sufficiently small $\mu > 0$ such that there exists $\lambda > 0$ such 
that~(\ref{formula for fast convergence}) holds for all sufficiently large $n$
and we let $\pi$ denote an arbitrary $\mu n$-rich partition $\pi$ of $[n]$ into 
parts $V_1, \ldots, V_l$.
Recall the definition of $\mbP_{n,\pi}(l,d)$ from
Definition~\ref{definition of richness}~(iii).
First we show that there are $q \in \mbbN$ and 
a first-order formula $\xi(x,y)$ in the language of graphs
such that if $\mcG \in \mbP_{n,\pi}(l,d)$ has the $q$-extension property, then
$\mcG \models \xi(v,w)$ if and only if $v$ and $w$ belong to the same part of the 
partition $\pi$.
This result is then used to show that the proportion of $\mcG \in \mbP_n(l,d)$ which
have a unique decomposition approaches 1 as $n \to \infty$
(Theorem~\ref{unique decompositions}).

\begin{lem}\label{existence of neighbours in another class}
Let $m_1, m_2 \in \mbbN$.
For all large enough $n$, if $\mcG = (V, E^{\mcG}) \in \mbP_{n,\pi}(l,d)$ 
has the $(m_1 + m_2)$-extension property with respect to $\pi$, then $\mcG$ has the following property:
\begin{itemize}
\item[] Whenever $p \in [l]$, $X \subseteq V \setminus V_p$, $|X| \leq m_1$,
then there are (at least) $m_2$
distinct vertices $v_1, \ldots, v_{m_2} \in V_p$ such that $v_i$ is adjacent to every member in $X$,
for $i = 1, \ldots, m_2$.
\end{itemize}
\end{lem}

\noindent
{\bf Proof.}
Let $m_1, m_2 \in \mbbN$.
Let $k = m_1 + m_2$ and suppose that $\mcG = (V, E^{\mcG}) \in \mbP_{n,\pi}(l,d)$ 
has the $k$-extension property.
Let $p \in [l]$ and let $X \subseteq V \setminus V_p$ satisfy $|X| \leq m_1$.
Now we define a suitable $\mcH = (X_1 \cup X_2 \cup Y, E^{\mcH})$ satisfying 
Assumption~\ref{assumption about H}, and then apply 
Lemma~\ref{second extension lemma}.

Let $X_1 = \es$, $X_2 = X$ and $Y = \{a_1, \ldots, a_{m_2}\}$
where $a_1, \ldots, a_{m_2}$ are new vertices.
Then let $X_1 \cup X_2 \cup Y$ be the vertex set of $\mcH$ and define the edge relation $E^{\mcH}$
as follows: $\mcH[X_2] = \mcG[X_2]$ 
(recall that $X_2 = X \subseteq V \setminus V_p$), $a_i \sim_{\mcH} w$
for every $i$ and every $w \in X_2$, 
and $a_i \not\sim_{\mcH} a_j$ if $i \neq j$.
Let $h_0$ denote the identity function on $X_2 = X_1 \cup X_2$ (recall that $X_1 = \es$), 
so $h_0$ is a strong embedding of $\mcH[X_1 \cup X_2]$ into $\mcG$ such that
$h_0(X_1) = \es \subseteq V_p$ and $h_0(X_2) = X_2 = X \subseteq V \setminus V_p$.
Moreover, as no vertex of $\mcG[V_p]$ has degree more than $d$ and
$Y$ is an independent set of cardinality $m_2$ it follows,
for large enough $n$, that $\mcG$ has at least $n^{1/4}$ different induced
subgraphs that are isomorphic to $\mcH[Y]$.
As $\mcG$ has the $k$-extension property and
$|X_1| + |X_2| + |Y| \leq m_1 + m_2 = k$, it follows that $h_0$ can be extended to
a strong embedding $h : \mcH \to \mcG$ such that $h(Y) \subseteq V_p$.
If $v_i = h(a_i)$ for $i = 1, \ldots, m_2$, then 
$v_i \sim_{\mcG} w$ for every $i$ and every $w \in h(X_2) = X_2 = X$.
\hfill $\square$
\\

\noindent
We will use the following:

\begin{observation}\label{observation about the number of neighbours in the same part}{\rm
Let $v, w_1, \ldots, w_s$ be distinct vertices of a graph $\mcG \in \mbP_{n,\pi}(l,d)$.
If all $w_1, \ldots, w_s$ are neighbours of $v$,
then at least $s - d$ of the vertices $w_1, \ldots, w_s$ 
do {\em not} belong to the same $\pi$-part as $v$.
}\end{observation}

\begin{defin}\label{definition of xi}{\rm
Let $m = (l+1)d + 1$ 
and let $\xi(x,y)$ denote the formula
\begin{align*}
\exists z_1 \ldots z_{(l-1)m} \Bigg[ \bigwedge_{1 \leq i < j \leq (l-1)m} z_i \neq z_j
\ \ &\wedge \ 
\bigwedge_{1 \leq i \leq (l-1)m} \Big( E(x,z_i) \wedge E(y,z_i) \Big) \\
&\wedge \ \
\bigwedge_{2 \leq k \leq l-1} \ \ \bigwedge_{i \leq (k-1)m < j} E(z_i,z_j) \Bigg].
\end{align*}
}\end{defin}

\begin{lem}\label{definability of the partition}
Let $m = (l+1)d + 1$ (as in Definition~\ref{definition of xi}) and $q = 2 + lm$.
If $\mcG \in \mbP_{n,\pi}(l,d)$ has the $q$-extension property with respect to $\pi$,
then, for all vertices $v$ and $w$ of $\mcG$,
$$\text{$v$ and $w$ belong to the same $\pi$-part} \ \  \Longleftrightarrow \ \ 
\mcG \models \xi(v,w).$$
\end{lem}

\noindent
{\bf Proof.}
Suppose that $\mcG = (V, E^{\mcG}) \in \mbP_{n,\pi}(l,d)$ has the 
$q$-extension property with respect to $\pi$, where $q = 2 + lm$
and $m = (l+1)d + 1$.
Also assume that $v$ and $w$ are vertices of $\mcG$.

First suppose that $v$ and $w$ belong to the same $\pi$-part of $V = [n]$.
For the sake of simplicity of notation, and without loss of generality,
suppose that $v, w \in V_1$.
In order to show that $\mcG \models \xi(v,w)$ we need to find
distinct vertices $u_1, \ldots, u_{(l-1)m}$ such that
\begin{align}\label{quantifiers in xi-0 replaced by parameters}
\mcG \models \ \bigwedge_{1 \leq i < j \leq (l-1)m} u_i \neq u_j \ &\wedge \ 
\bigwedge_{1 \leq i \leq (l-1)m} \Big( E(v,u_i) \ \wedge \ E(w,u_i) \Big) \\ 
& \wedge \
\bigwedge_{2 \leq k \leq l-1} \ \ \bigwedge_{i \leq (k-1)m < j} E(u_i,u_j). \nonumber
\end{align}
This can be proved by showing, by induction, that for every $t = 1, \ldots, l-1$,
there are vertices $u_1, \ldots, u_{tm}$ such that
\begin{equation}\label{all u belong to the same pi-part}
\text{for every $k = 1, \ldots, t$, we have $u_{(k-1)m +1}, \ldots, u_{km} \in V_{k+1}$, and} 
\end{equation}
\begin{align}\label{the other conditions on the u}
\mcG \models \ 
\bigwedge_{1 \leq i < j \leq tm} u_i \neq u_j \ &\wedge \
\bigwedge_{1 \leq i \leq tm} \Big( E(v,u_i) \ \wedge \ E(w,u_i) \Big) \\ 
&\wedge \
\bigwedge_{2 \leq k \leq t} \ \ \bigwedge_{1 \leq i \leq (k-1)m < j \leq tm} E(u_i,u_j). \nonumber
\end{align}
In the base case $t=1$, we need to find $u_1, \ldots, u_m$ such
that~(\ref{all u belong to the same pi-part})
and~(\ref{the other conditions on the u}) hold for $t=1$.
But the existence of such $u_1, \ldots, u_m \in V_2$ is guaranteed by 
Lemma~\ref{existence of neighbours in another class}
and the assumption that $\mcG$ has the $q$-extension property.

In the inductive step we assume that $t < l-1$ and that there are
$u_1, \ldots, u_{tm}$ such that~(\ref{all u belong to the same pi-part})
and~(\ref{the other conditions on the u}) hold.
Recall the assumption that $v,w \in V_1$. 
Again, Lemma~\ref{existence of neighbours in another class}
and the assumption that $\mcG$ has the $q$-extension property
implies that there are
$u_{tm+1}, \ldots, u_{(t+1)m} \in V_{t+2}$ such
that~(\ref{all u belong to the same pi-part})
and~(\ref{the other conditions on the u}) hold if $t$ is replaced by $t+1$.

It remains to prove that if $\mcG \models \xi(v,w)$,
then $v$ and $w$ belong to the same $\pi$-part.
So suppose that $\mcG \models \xi(v,w)$, which implies that there
are distinct vertices $u_1, \ldots, u_{(l-1)m}$
such that~(\ref{quantifiers in xi-0 replaced by parameters}) holds.
For a contradiction, suppose that $v$ and $w$ do not belong to the same $\pi$-part.
By Observation~\ref{observation about the number of neighbours in the same part}
and the choice of $m = (l+1)d + 1$, 
it follows that there are $i_1, \ldots, i_{l-1}$ such that, for every $k = 1, \ldots, l-1$
(recall that $l \geq 2$),
$(k-1)m < i_k \leq km$ and $u_{i_k}$
does not belong to the same $\pi$-part as any of $v, w, u_{i_1}, \ldots, u_{i_{k-1}}$.
As $v$ and $w$ do not belong to the same $\pi$-part (by assumption) this contradicts
that there are only $l$ $\pi$-parts.
\hfill $\square$

\begin{defin}\label{definition of definable partition}{\rm
Let $\mcG = (V, E^{\mcG}) \in \mbP(l,d)$.\\
(i) A relation $R \subseteq V^k$ is called {\em definable in $\mcG$} by a formula
$\varphi(x_1, \ldots, x_k)$ if for all $(v_1, \ldots, v_k) \in V^k$,
$(v_1, \ldots, v_k) \in R$ $\Longleftrightarrow$ $\mcG \models \varphi(v_1, \ldots, v_k)$.\\
(ii) A partition $\pi$ of $V$ is called {\em definable in $\mcG$} by a first-order formula
$\varphi(x_1,x_2)$ if the
relation `$v$ and $w$ belong to the same $\pi$-part' is definable by $\varphi(x_1,x_2)$. 
}\end{defin}

\begin{theor}\label{unique decompositions}
For every $k \in \mbbN$, the proportion of $\mcG \in \mbP_n(l,d)$
with the following properties approaches 1 as $n \to \infty$:
\begin{itemize}
	\item[(i)] $\mcG$ has a unique decomposition,
	\item[(ii)] The formula $\xi(x,y)$
	(from Definition~\ref{definition of xi}) 
	defines the partition on which the unique decomposition of $\mcG$ is based, 
	and this partition
	is $\mu n$-rich.
	\item[(iii)] $\mcG$ has the $k$-extension property with respect to the partition on
	which its unique decomposition is based.
\end{itemize}
\end{theor}

\noindent
{\bf Proof.}
It suffices to prove the proposition for all sufficiently large $k$.
So we assume that $k \geq 2 + lm$, where $m = (l+1)d + 1$ (which will allow us to 
use Lemma~\ref{definability of the partition}).
For every $\mu n$-rich partition $\pi$ of $[n]$ let $\mbX_{n,\pi}$ denote the set of
$\mcG \in \mbP_{n,\pi}(l,d)$ which have the $k$-extension property with respect to $\pi$.
By Corollary~\ref{corollary to second extension property},
there is $\varepsilon_k : \mbbN \to \mbbR$ such that $\lim_{n\to\infty} \varepsilon_k(n) = 0$
and 
\begin{align}
&\text{for every $n$ and $\mu n$-rich partition $\pi$ of $[n]$ into $l$ parts,} \nonumber \\
&|\mbX_{n,\pi}| \big/ |\mbP_{n,\pi}(l,d)| \geq 1 - \varepsilon_k(n).
\label{lower bound on number of graphs with k-extension property}
\end{align}

\noindent
{\em Claim.} If $\pi_1$ and $\pi_2$ are different $\mu n$-rich partitions of $[n]$
into $l$ parts, then \\$\mbX_{n,\pi_1} \cap \mbX_{n,\pi_2} = \es$.
\\

\noindent
{\em Proof of Claim.}
Suppose that $\pi_1$ and $\pi_2$ are $\mu n$-rich partitions of $[n]$
into $l$ parts and that $\mcG \in \mbX_{n,\pi_1} \cap \mbX_{n,\pi_2}$.
Lemma~\ref{definability of the partition} and the choice of $k$ (being large enough)
implies that, for all $v,w \in [n]$,
\begin{align*}
&\text{$v$ and $w$ belong to the same part with respect to $\pi_1$}\\
\Longleftrightarrow \
&\mcG \models \xi(v,w)\\
\Longleftrightarrow \
&\text{$v$ and $w$ belong to the same part with respect to $\pi_2$.}
\end{align*}
Hence, $\pi_1 = \pi_2$.
\hfill $\square$
\\

\noindent
Let 
$$\mbX_n \ = \ \bigcup_{\pi \text{ $\mu n$-rich}} \mbX_{n,\pi},$$
where the union ranges over all $\mu n$-rich partitions $\pi$ of $[n]$, and let
$\mbP_n^*$ be the set of $\mcG \in \mbP_n(l,d)$ such that {\em every} decomposition of $\mcG$ 
is based on a partition of $[n]$ that is $\mu n$-rich.
By the claim and~(\ref{lower bound on number of graphs with k-extension property}) we get
\begin{align}\label{lower bound on number of graphs with k-extension property 2}
|\mbX_n| \ &= \ \sum_{\pi \text{ $\mu n$-rich}} |\mbX_{n,\pi}| \ \geq \ 
\sum_{\pi \text{ $\mu n$-rich}} (1 - \varepsilon_k(n))|\mbP_{n,\pi}|  \\
&= \ 
(1 - \varepsilon_k(n)) \sum_{\pi \text{ $\mu n$-rich}} |\mbP_{n,\pi}| \ \geq \ 
(1 - \varepsilon_k(n))|\mbP_n^*|, \nonumber
\end{align}
where the sums range over all $\mu n$-rich partitions $\pi$ of $[n]$ into $l$ parts.
By the choice of $\mu$ and Corollary~\ref{almost all graphs have a partition with large parts},
\begin{equation}\label{number of graphs with only mu n-rich partitions approaches one}
\lim_{n\to\infty} \frac{|\mbP_n^*|}{|\mbP_n(l,d)|} \ = \ 1.
\end{equation}
We now get
$$1 \ \geq \ 
\frac{|\mbX_n|}{|\mbP_n(l,d)|} \ = \ \frac{|\mbX_n|}{|\mbP_n^*|} \cdot \frac{|\mbP_n^*|}{|\mbP_n(l,d)|} \ 
\overset{(\ref{lower bound on number of graphs with k-extension property 2})}{\geq} \ 
(1 - \varepsilon_k(n))\frac{|\mbP_n^*|}{|\mbP_n(l,d)|},$$
so by~(\ref{number of graphs with only mu n-rich partitions approaches one}),
$$\lim_{n\to\infty}\frac{|\mbX_n|}{|\mbP_n(l,d)|} \ = \ 1,$$
and together with~(\ref{number of graphs with only mu n-rich partitions approaches one}),
this implies
$$\lim_{n\to\infty}\frac{|\mbX_n \cap \mbP_n^*|}{|\mbP_n(l,d)|} \ = \ 1.$$
Thus it suffices to prove that every $\mcG \in \mbX_n \cap \mbP_n^*$
satsifies~(i)--(iii) of the proposition.
So let $\mcG \in  \mbX_n \cap \mbP_n^*$.
As $\mcG \in \mbP_n^*$, every decomposition of $\mcG$ is based on a $\mu n$-rich partition of $[n]$
into $l$ parts. Let $\pi_1$ and $\pi_2$ be two $\mu n$-rich partitions such that
$\mcG$ has decompositions based on $\pi_1$ and on $\pi_2$.
Then $\mcG \in \mbX_{n,\pi} \cap \mbX_{n,\pi}$ so by the claim, $\pi_1 = \pi_2$.
Hence all decompositions of $\mcG$ are based on the same partition, and thus,
by Observation~\ref{observation about decompositions}~(c), $\mcG$ has a unique decomposition.
As $\mcG \in \mbP_n^*$ this partition, say $\pi$, is $\mu n$-rich, and as $\mcG \in \mbX_{n,\pi}$,
it follows from Lemma~\ref{definability of the partition} and the choice of $k$
that $\xi(x,y)$ defines $\pi$.
\hfill $\square$

\section{A limit law}\label{a limit law}

\noindent
In this section we prove the main result of this article, Theorem~\ref{limit theorem},
in a slightly different formulation compared with its statement in Section~\ref{introduction}.
In Section~\ref{forbidden subgraphs} we use it to get
limit laws for $\mcH$-free graphs for certain types of $\mcH$.
\medskip

\noindent
{\bf Theorem~\ref{limit theorem}.} 
{\em
Let $l \geq 1$ be an integer.\\
(i) $\mbP(l,1)$ has a zero-one law.\\
(ii) For every $d \geq 2$, $\mbP(l,d)$ has a limit law, but not zero-one law.
}
\medskip

\noindent
We first prove part~(ii) of Theorem~\ref{limit theorem}
and then, in Section~\ref{proof of part one of limit law for almost l-partite graphs},
sketch the much easier proof of part one.

\subsection{Proof of part (ii) of Theorem~\ref{limit theorem}}
\label{proof of part two of limit law for almost l-partite graphs}

Suppose that $d \geq 2$. As the case $l = 1$ is proved in \cite{Kop12} we also
assume that $l \geq 2$.
Let $\varphi$ be an arbitrary first-order sentence in the language of graphs.
We need to prove that the quotient
$$\frac{|\{\mcG \in \mbP_n(l,d) : \mcG \models \varphi \}|}{|\mbP_n(l,d)|}$$
converges as $n \to \infty$.
Suppose that the quantifier rank of $\varphi$ is at most $k$, where $k \geq 1$.
(See for example \cite{EF} for the definition of quantifier rank.)

\begin{defin}\label{definition of k-elementary equivalence}{\rm
For graphs $\mcG_1$ and $\mcG_2$ let $\mcG_1 \equiv_k \mcG_2$ mean 
that $\mcG_1$ and $\mcG_2$ satisfy exactly the same first-order sentences
with quantifier rank at most $k$.
}\end{defin}

\noindent
Note that `$\equiv_k$' is an equivalence relation on $\mbP(l,d) = \bigcup_{n \in \mbbN^+} \mbP_n(l,d)$.
As the language of graphs has a finite relational vocabulary it follows that
$\equiv_k$ has only finitely many equivalence classes (e.g. \cite{EF}).
Therefore, to prove that $\mbP(l,d)$ has a limit law, 
it suffices to prove that for every $\equiv_k$-class $\mbE$, the quotient
$$\frac{|\mbE \cap \mbP_n(l,d)|}{|\mbP_n(l,d)|}$$
converges as $n \to \infty$.
This is done in Corollary~\ref{convergence of any subset of classes},
which finishes the proof of the first claim of 
part~(ii) of Theorem~\ref{limit theorem}.
From the proof it is easy to deduce that a zero-one law does not hold, which is explained
after the proof of Corollary~\ref{convergence of any subset of classes}.

Recall that $\mbP_n(1,d)$ (the case $l = 1$) denotes the set of graphs with vertex set $[n]$
in which every vertex has degree at most $d$.
We will use results from \cite{Kop12} about the asymptotic structure of graphs in $\mbP_n(1,d)$,
stated as Theorems~\ref{properties that are almost surely true}
and~\ref{distribution of poisson objects} below.

\begin{theor}\label{properties that are almost surely true} \cite{Kop12}
Suppose that $d \geq 2$ and $s, t > 0$ are integers and that $0 < \varepsilon < d$.
The proportion of $\mcG \in \mbP_n(1,d)$ with properties~(1)--(6) below
approaches 1 as $n \to \infty$. If $d \geq 3$ then the
proportion of $\mcG \in \mbP_n(1,d)$ with properties~(1)--(7) approaches 1 as $n \to \infty$.
\begin{enumerate}
\item There is no vertex with degree degree less than $d-2$.

\item There are between $\sqrt{(d - \varepsilon)n}$ and 
$\sqrt{(d + \varepsilon)n}$ vertices with degree $d-1$.

\item If $p, q \leq s$ then there are no $p$-cycle and different $q$-cycle within
distance at most $t$ of each other.

\item If $p \leq s$ then there are no vertex $v$ with degree less than $d$ and 
$p$-cycle within distance at most $t$ of each other.
In particular, no $p$-cycle contains
a vertex of degree less than $d$.

\item There do not exist distinct vertices $v_1, v_2, v_3$ all of
which have degree at most $d-1$ such that for all distinct $i,j \in \{1,2,3\}$,
$\dist_{\mcG}(v_i,v_j) \leq t$.

\item There do not exist distinct vertices $v$ and $w$ such that $\deg_{\mcG}(v) \leq d-1$,
$\deg_{\mcG}(w) \leq d-2$ and $\dist_{\mcG}(v,w) \leq t$.

\item Every connected component has more than $t$ vertices.
\end{enumerate}
\end{theor}

\begin{defin}\label{definition of union of graphs}{\rm
For graphs $\mcG_i = (V_i, E^{\mcG_i})$, $i = 1, \ldots, m$,
$\bigcup_{i=1}^m \mcG_i$ denotes the graph $\mcG = (V,E^{\mcG})$
where $V = \bigcup_{i=1}^m V_i$ and $E^{\mcG} = \bigcup_{i=1}^m E^{\mcG_i}$.
}\end{defin}

\begin{rem}\label{remark on construction of the graphs and probabilities}{\rm
Let $\pi$ denote the partition $V_1, \ldots, V_l$ of $V = [n]$.
By definition of $\mbP_{n,\pi}(l,d)$ (Definition~\ref{definition of richness}~(iii)) and
Observation~\ref{observation about decompositions}~(c),
every $\mcG \in \mbP_{n,\pi}(l,d)$ can be constructed in a 
unique way by choosing $\mcG_1 \in \mbP_{n,\pi}(l,0)$ and then,
for $i = 1, \ldots, l$, choosing $\mcG'_i = (V_i, E^{\mcG'_i})$ 
in which every vertex has degree at most $d$ (and every choice is independent of the 
previous choices) and letting 
$$\mcG = \mcG_1 \cup \bigcup_{i=1}^l \mcG'_i.$$
Conversely, every graph that is constructed by this procedure belongs to $\mbP_{n,\pi}(l,d)$.
Therefore,
\begin{equation}\label{number of elements in P-n-pi as a product}
|\mbP_{n,\pi}(l,d)| \ = \ |\mbP_{|V_1|}(1,d)| \cdots |\mbP_{|V_l|}(1,d)| \cdot |\mbP_{n,\pi}(l,0)|,
\end{equation}
}\end{rem}

\noindent
The set defined in the next definition contains the typical (for large enough $n$)
graphs in $\mbP_n(l,d)$.

\begin{defin}\label{definition of P-t}{\rm
{\em Fix some $0 < \varepsilon < d$ for the rest of this section.
Also fix some sufficiently small $\mu > 0$ such that there exists $\lambda > 0$ such 
that~(\ref{formula for fast convergence}) of Section~\ref{introduction} 
holds for all sufficiently large $n$.}
Let $\mbP^k_n(l,d)$ be the set of all $\mcG \in \mbP_n(l,d)$ such that
\begin{itemize}
\item[(i)] \ $\mcG$ has a decomposition which is based on a $\mu n$-rich partition
$V_1, \ldots, V_l$ of $[n]$ and this partition is defined by $\xi(x,y)$,
\item[(ii)] \ for $s = 5^k$, $t = 5^{k+1}$ and every $i \in [l]$ properties (1)--(6) of
Theorem~\ref{properties that are almost surely true} hold for $\mcG[V_i]$, and
\item[(iii)] \ if $d \geq 3$ then, for $t = 5^{k+1}$, also property~(7)
of Theorem~\ref{properties that are almost surely true} holds for $\mcG[V_i]$ for every $i \in [l]$.
\end{itemize}
Let $\mbP^k(l,d) = \bigcup_{n = 1}^{\infty}\mbP^k_n(l,d)$.
}\end{defin}

\begin{lem}\label{P-t is almost everything}
$\displaystyle \lim_{n\to\infty} |\mbP^k_n(l,d)| \big/ |\mbP_n(l,d)| = 1$. 
\end{lem}

\noindent
{\bf Proof.}
By the choice of $\mu$ in Definition~\ref{definition of P-t}
and~Theorem~\ref{unique decompositions}, the proportion of $\mcG \in \mbP_n(l,d)$
for which~(i) of Definition~\ref{definition of P-t} holds
approaches~1 as $n \to \infty$.
Now consider any $\mu n$-rich partition $\pi$ of $V = [n]$ with parts $V_1, \ldots, V_l$,
so $|V_i| \geq \mu n$ for all $i \in [l]$.
From~(\ref{number of elements in P-n-pi as a product}) 
and 
Theorem~\ref{properties that are almost surely true} 
it follows that the proportion of $\mcG \in \mbP_{n,\pi}(l,d)$ such 
that for every $i \in [l]$ properties (1)--(6) (and~(7) if $d \geq 3$) of
Theorem~\ref{properties that are almost surely true} 
hold for $\mcG[V_i]$ approaches 1 as $n \to \infty$.
Since $|V_i| \geq \mu n$ for all $i \in [l]$, the rate of convergence
depends only on $l$ and $d$, as $\mu$ depends only on $l$.
Therefore $|\mbP^k_n(l,d)| \big/ |\mbP_n(l,d)| \to 1$ as $n \to \infty$.
\hfill $\square$
\\

\noindent
Now we define an equivalence relation on $\mbP(l,d)$ which distinguishes
whether a graph $\mcG$ belongs to $\mbP^k(l,d)$ or not, and
if the answer is `yes' then it distinguishes the number of
vertices of degree $d-2$, the number of $i$-cycles for $i \leq 5^{k}$
and the number of $i$-paths with endpoints of degree $d-1$ for $i \leq 5^{k}$
in $\mcG[V_j]$, for each part $V_j$ of the partition on which the unique
decomposition of $\mcG$ is based.

\begin{defin}\label{definition of poisson equivalence}{\rm
We define an equivalence relation `$\approx_k$' on $\mbP(l,d)$ as follows:\\
$\mcG \approx_k \mcH$ if and only if
\begin{itemize}
\item[\textbullet] \  either $\mcG, \mcH \notin \mbP^k(l,d)$ or
\item[\textbullet] \ $\mcG, \mcH \in \mbP^k(l,d)$ and if 
$V_1, \ldots, V_l$ is the partition of the vertex set of $\mcG$ defined by $\xi(x,y)$ on which some 
decomposition of $\mcG$ is based and 
$W_1, \ldots, W_l$ is the partition of the vertex set of $\mcH$ defined by $\xi(x,y)$ on which some 
decomposition of $\mcH$ is based,
then there is a permutation $\sigma$ of $[l]$ such that, for every $i \in [l]$,
\begin{itemize}
\item[(a)] \ $\mcG[V_i]$ and $\mcH[W_{\sigma(i)}]$ have the same number of vertices with degree $d-2$,
\item[(b)] \ for every $j = 3, \ldots, 5^{k}$, $\mcG[V_i]$ and $\mcH[W_{\sigma(i)}]$ have the same number of $j$-cycles, and 
\item[(c)] \ for every $j = 1, \ldots, 5^{k}$, $\mcG[V_i]$ and $\mcH[W_{\sigma(i)}]$ have the same number of $j$-paths
with both endpoints of degree $d-1$.
\end{itemize}
\end{itemize}
}\end{defin}

\noindent
The next result from \cite{Kop12} will be used to show that
for every $\approx_k$-equivalence class $\mbC$, the quotient
$|\mbC \cap \mbP_n(l,d)| \big/ |\mbP_n(l,d)|$ converges when $n \to \infty$.

\begin{theor}\label{distribution of poisson objects} \cite{Kop12}
Let $t \geq 3$ be an integer. 
There are positive $\lambda_3, \ldots, \lambda_t, \mu_1, \ldots, \mu_t \in \mbbQ$ 
such that for all $q, r_3, \ldots, r_t$, $s_1, \ldots, s_t \in \mbbN$ 
the proportion of $\mcG \in \mbP_n(1,d)$ such that 
\begin{itemize}
\item[(a)] $\mcG$ has exactly $q$ vertices with degree $d-2$,

\item[(b)] for $i = 3, \ldots, t$, $\mcG$ has exactly $r_i$ $i$-cycles,  and 

\item[(c)] for $i = 1, \ldots, t$, $\mcG$ has exactly $s_i$ $i$-paths with both endpoints
of degree $d-1$
\end{itemize}
approaches
$$\frac{(d-1)^q \ e^{-(d-1)}}{q!} \ 
\Bigg(\prod_{i = 3}^t \frac{(\lambda_i)^{r_i} \ e^{-\lambda_i}}{r_i!}\Bigg)
\Bigg(\prod_{i = 1}^t \frac{(\mu_i)^{s_i} \ e^{-\mu_i}}{s_i!}\Bigg)
\quad \text{ as } \ n \to \infty.$$
\end{theor}

\noindent
In other words, Theorem~\ref{distribution of poisson objects}
says that the random variables which, for a random $\mcG \in \mbP_n(1,d)$, count
the number of vertices with degree $d-2$,
the number of $i$-cycles for $3 \leq i \leq t$ and
the number of $i$-paths with both endpoints of degree $d-1$ for $1 \leq i \leq t$
have independent Poisson distributions, asymptotically.
Note that since the Poisson distribution is a probability distribution it follows
that the sum of all numbers as in the conclusion of
Theorem~\ref{distribution of poisson objects} when 
$q, r_3, \ldots, r_t$, $s_1, \ldots, s_t$ ranges over all natural numbers is 1.

\begin{lem}\label{convergence of poisson classes}
For every equivalence class $\mbC$ of the relation `$\approx_k$' there is 
a constant $0 \leq c(\mbC) \leq 1$ such that 
$$\lim_{n\to\infty} \frac{|\mbC \cap \mbP_n(l,d)|}{|\mbP_n(l,d)|} \ = \ c(\mbC).$$
Moreover, $c(\mbC) = 0$ if and only if $\mbC = \mbP(l,d) \setminus \mbP^k(l,d)$.
If $(\mbC_i : i \in \mbbN)$ is an enumeration of all $\approx_k$-classes then
$\displaystyle \sum_{i=1}^{\infty} c(\mbC_i) = 1$.
\end{lem}

\noindent
{\bf Proof.}
If $\mbC = \mbP(l,d) \setminus \mbP^k(l,d)$
then the conclusion holds with $c = 0$, by Lemma~\ref{P-t is almost everything}.
If $\mbC \neq \mbP(l,d) \setminus \mbP^k(l,d)$ then $\mbC \subseteq \mbP^k(l,d)$
and the conclusion holds because of~(\ref{number of elements in P-n-pi as a product}),
Lemma~\ref{P-t is almost everything} and
Theorem~\ref{distribution of poisson objects}.
\hfill $\square$

\begin{defin}\label{definition of poisson object}{\rm
If $\mcG$ is a graph in which every vertex has degree at
most $d$, then let a {\em small Poisson object} of $\mcG$ denote any one of
\begin{itemize}
\item[(a)] \ a vertex with degree $d-2$, or
\item[(b)] \ an $i$-cycle where $3 \leq i \leq 5^{k}$, or
\item[(c)] \ an $i$-path with both endpoints of degree $d-1$ where $1 \leq i \leq 5^{k}$.
\end{itemize}
}\end{defin}

\begin{defin}\label{definition of neighbourhood}{\rm
Let $\mcG = (V, E^{\mcG})$, $A \subseteq V$ and $t \in \mbbN$.\\
(i) $N_{\mcG}(A,t) = \{ v \in V : \ \dist_{\mcG}(A, v) \leq t\}$.
Note that $N_{\mcG}(A,0) = A$.\\
(ii) $\mcN_{\mcG}(A,t) = \mcG\big[N_{\mcG}(A,t)\big]$.
Note that $\mcN_{\mcG}(A,0) = \mcG[A]$.\\
(iii) If every vertex in $\mcG$ has degree at most $d$,
then let $NP(\mcG, t)$ be the set of vertices $v$ of $\mcG$ such that
the distance from $v$ to a small Poisson object (of $\mcG$) is at most $t$.\\
(iv) If every vertex in $\mcG$ has degree at most $d$,
then let $\mc{NP}(\mcG, t) = \mcG[NP(\mcG, t)]$.
}\end{defin}

\begin{rem}\label{remark on definition of poisson equivalence}{\rm
Recall Definition~\ref{definition of union of graphs} and
observe that if $\mcG, \mcH \in \mbP^k(l,d)$, then $\mcG \approx_k \mcH$
if and only if 
there is an isomorphism
$$f : \bigcup_{i=1}^l \mc{NP}(\mcG[V_i], 5^{k}) \ \to \ \bigcup_{i=1}^l \mc{NP}(\mcH[W_i], 5^{k})$$
where $V_1, \ldots, V_l$ and $W_1, \ldots, W_l$ denote the $\xi(x,y)$-classes of $\mcG$ and $\mcH$,
respectively, and $f$ preserves $\xi(x,y)$ in the sense that whenever $v$ and $w$ 
are in the domain of $f$, then $\mcG \models \xi(v,w)$ if and only if $\mcH \models \xi(f(v),f(w))$.
}\end{rem}

\noindent
Note that if $\mcG \in \mbP^k(l,d)$ and the $V_1, \ldots, V_l$ is the partition
defined by $\xi(x,y)$, then the graph
$$\mcG\bigg[\bigcup_{i=1}^l NP(\mcG[V_i], 5^{k})\bigg]$$
is the result of adding to $\bigcup_{i=1}^l \mc{NP}(\mcG[V_i], 5^{k})$ 
all edges of $\mcG$ that connect
$NP(\mcG[V_i], 5^{k})$ and $NP(\mcG[V_j], 5^{k})$ for all distinct $i,j \in [l]$.

\begin{defin}\label{definition of strong poisson equivalence}{\rm
We define an equivalence relation `$\approx_k^+$' on $\mbP(l,d)$ as follows:\\
$\mcG \approx_k^+ \mcH$ if and only if
\begin{itemize}
\item[\textbullet] \ $\mcG \approx_k \mcH$ and 

\item[\textbullet] \ if $\mcG, \mcH \in \mbP^k(l,d)$ then there is an isomorphism 
$$f : \mcG\bigg[\bigcup_{i=1}^l NP(\mcG[V_i], 5^{k})\bigg] \ \to \ 
\mcH\bigg[\bigcup_{i=1}^l NP(\mcH[W_i], 5^{k})\bigg]$$
where $V_1, \ldots, V_l$ and $W_1, \ldots, W_l$ denote the $\xi(x,y)$-classes of $\mcG$ and $\mcH$,
respectively, and $f$ preserves $\xi(x,y)$ in the sense that whenever $v$ and $w$ 
are in the domain of $f$, then $\mcG \models \xi(v,w)$ if and only if $\mcH \models \xi(f(v),f(w))$.
\end{itemize}
}\end{defin}

\noindent
Note that the equivalence relation $\approx_k^+$ refines $\approx_k$, that is,
every $\approx_k^+$-class is included in some $\approx_k$-class.
Also note that every $\approx_k$-class is divided into finitely many $\approx_k^+$-classes,
and that there are an infinite (but countable) number of $\approx_k$-classes,
and hence an infinite (but countable) number of $\approx_k^+$-classes.

\begin{lem}\label{convergence of strong poisson classes modulo poisson classes}
Let $\mbC$ be an $\approx_k$-class such that $\mbC \subseteq \mbP^k(l,d)$
and let $\mbD$ be an $\approx_k^+$-class such that $\mbD \subseteq \mbC$.
Then there is a constant $c(\mbD, \mbC) > 0$, depending only on $l$, $d$, $k$ and $\mbD$,
such that 
$$\lim_{n \to \infty} \frac{|\mbD \cap \mbP_n(l,d)|}{|\mbC \cap \mbP_n(l,d)|} \ = \  c(\mbD, \mbC).$$
In particular, if $\mbD_1, \ldots, \mbD_m$ enumerates all $\approx_k^+$-classes
that are included in $\mbC$, then 
$$\sum_{i=1}^m c(\mbD_i, \mbC) = 1.$$
\end{lem}

\noindent
{\bf Proof.} 
Suppose that $\mbC$ is an equivalence class of `$\approx_k$' such that $\mbC \subseteq \mbP^k(l,d)$
and let $\mbD$ be an equivalence class of `$\approx_k^+$' such that $\mbD \subseteq \mbC$.
By Lemma~\ref{convergence of poisson classes}, there is $c > 0$ such that
$|\mbC \cap \mbP_n(l,d)| \big/ |\mbP_n(l,d)| \to c$ as $n \to \infty$.
For $\mcG \in \mbD \cap \mbP_n(l,d)$, let $p$ be the number of ways in which edges can be
added to $\bigcup_{i=1}^l \mc{NP}(\mcG[V_i], 5^k)$, where $V_1, \ldots, V_l$ are the
$\xi$-classes of $\mcG$, in such a way that
the resulting graph is isomorphic, via an isomorphism preserving the partition $V_1, \ldots, V_l$,
to $\mcG\Big[\bigcup_{i=1}^l NP(\mcG[V_i], 5^k)\Big]$.
Note that $p$ depends only on $\mbD$ (and not on $n$ or the particular graph $\mcG$ from $\mbD$).
For $\mcG \in \mbC \cap \mbP_n(l,d)$ let $q$ be the total number of ways in which edges can be added
between $NP(\mcG[V_i], 5^k)$ and $NP(\mcG[V_j], 5^k)$, where the $V_1, \ldots, V_l$ 
are the $\xi$-classes of $\mcG$, for all possible distinct $i, j \in [l]$.
Also let $r = \Big|\bigcup_{i=1}^l NP(\mcG[V_i], 5^k)\Big|$.
Then $q$ and $r$ depend only on $\mbC$.
We show that
\begin{equation}\label{D divided by C}
\lim_{n\to\infty} \ \frac{\big|\mbD \cap \mbP_n(l,d)\big|}{\big|\mbC \cap \mbP_n(l,d)\big|} \ = \ \frac{p}{q}.
\end{equation}
For $\mcG \in \mbC \cap \mbP_n(l,d)$ with $\xi$-classes $V_1, \ldots, V_l$
(which is a $\mu n$-rich partition of $[n]$),
let $[\mcG]_n$ be the set of $\mcH \in \mbC \cap \mbP_n(l,d)$
such that $\mcH$ has a decomposition based on $V_1, \ldots, V_l$,
$$\bigcup_{i=1}^l \mc{NP}(\mcH[V_i], 5^k) \ = \ 
\bigcup_{i=1}^l \mc{NP}(\mcG[V_i], 5^k)$$
and if at least one of $v$ or $w$ does {\em not} belong to $\bigcup_{i=1}^l NP(\mcG[V_i], 5^k)$, 
then $v \sim_{\mcG} w$ $\Longleftrightarrow$ $v \sim_{\mcH} w$.
Note that $\big|[\mcG]_n\big|$ has a finite bound depending only on $\mbC$.
By Lemma~\ref{definability of the partition}, there is $m$ such that
whenever $\mcG \in \mbP_n(l,d)$ has a decomposition based on a $\mu n$-rich partition $\pi$
and $\mcG$ has the $m$-extension property with respect to $\pi$, then $\xi(x,y)$ defines $\pi$.
By Lemma~\ref{removing some edges from a graph with large extension property},
if $\mcG \in \mbC \cap \mbP_n(l,d)$ and at least one member of $[\mcG]_n$ has the 
$(m + 2\binom{r}{2})$-extension property with respect to the partition $V_1, \ldots, V_l$
defined by $\xi(x,y)$,
then all members of $[\mcG]_n$ have the $m$-extension property with respect to $V_1, \ldots, V_l$.
So if $\mcG \in \mbC \cap \mbP_n(l,d)$ and at least one member of 
$[\mcG]_n$ has the $(m + 2\binom{r}{2})$-extension property
with respect to the partition defined by $\xi(x,y)$,
then the proportion of $\mcH \in [\mcG]_n$ which belong to $\mbD$ is exactly $p/q$.
By Theorem~\ref{unique decompositions} the proportion of $\mcG \in \mbP_n(l,d)$
in which $\xi(x,y)$ defines a partition (equivalence relation) and $\mcG$
has the $(m + 2\binom{r}{2})$-extension property with respect to this partition
approaches 1 as $n \to \infty$,
Since $|\mbC \cap \mbP_n(l,d)| \big/ |\mbP_n(l,d)| \to c > 0$ as $n \to \infty$,
it follows that the proportion of $\mcH \in \mbC \cap \mbP_n(l,d)$ which have
the $(m + 2\binom{r}{2})$-extension property approaches 1 as $n \to \infty$.
As mentioned above, there is a finite bound, say $\beta$, depending only on $\mbC$
such that for every $\mcG \in \mbC \cap \mbP_n(l,d)$, 
$\big|[\mcG]_n\big| \leq \beta$. 
It follows that the proportion of $\mcH \in \mbC \cap \mbP_n(l,d)$ which belong to
$\mbD$ approaches $p/q$ as $n \to \infty$, so~(\ref{D divided by C}) is proved.
\hfill $\square$

\begin{cor}\label{convergence of strong poisson classes}
Let $\mbD$ be an $\approx_k^+$-class.
Then there is a constant $c(\mbD)$, depending only on $l$, $d$, $k$ and $\mbD$,
such that 
$$\lim_{n\to\infty} \frac{|\mbD \cap \mbP_n(l,d)|}{|\mbP_n(l,d)|} \ = \ c(\mbD).$$
Moreover, $c(\mbD) = 0$ if and only if 
$\mbD = \mbP(l,d) \setminus \mbP^k(l,d)$.
If $(\mbD_i : i \in \mbbN)$ enumerates all $\approx_k^+$-classes,
then $\sum_{n=0}^{\infty} c(\mbD_i) = 1$.
\end{cor}

\noindent
{\bf Proof.}
Immediate from Lemmas~\ref{convergence of poisson classes}
and~\ref{convergence of strong poisson classes modulo poisson classes}.
\hfill $\square$
\\

\begin{lem}\label{convergence of k-equivalence class cut with strong poisson class}
Let $\mbD$ be any $\approx_k^+$-class such that 
$$\lim_{n\to\infty} \ \frac{|\mbD \cap \mbP_n(l,d)|}{|\mbP_n(l,d)|} \ = \ c(\mbD) > 0.$$
Then there is an $\equiv_k$-class $\mbE$ such that
$$\lim_{n\to\infty} \ \frac{|\mbE \cap \mbD \cap \mbP_n(l,d)|}{|\mbD \cap \mbP_n(l,d)|} \ = \ 1,$$
and consequently
$$\lim_{n\to\infty} \ \frac{|\mbE \cap \mbD \cap \mbP_n(l,d)|}{|\mbP_n(l,d)|} \ = \ c(\mbD).$$
\end{lem}

\noindent
{\bf Proof.}
Let $\mbD$ be any $\approx_k^+$-class such that 
$|\mbD \cap \mbP_n(l,d)| \big/ |\mbP_n(l,d)| \to c(\mbD) > 0$ as $n \to \infty$.
Then, by the definitions of $\approx_k^+$ and $\approx_k$, $\mbD \subseteq \mbP^k(l,d)$,
for every $i \in [l]$ there are 
$q_i, r_{i,3}, \ldots, r_{i,5^k}, s_{i,1}, \ldots, s_{i,5^k} \in \mbbN$
such that whenever $\mcG_1 = (V, E^{\mcG_1}) \in \mbD$ and
$\mcG_2 = (W, E^{\mcG_2}) \in \mbD$, then the following hold:
\begin{itemize}
\item[($*$)] The parts of the partitions of $V$ and of $W$ defined by $\xi$ in $\mcG_1$ and $\mcG_2$, respectively, 
can be ordered as $V_1, \ldots, V_l$ and as $W_1, \ldots, W_l$ 
in such a way that, for every $i \in [l]$,
\begin{itemize}
\item[(a)] both $\mcG_1[V_i]$ and $\mcG_2[W_i]$ have exactly $q_i$ vertices with degree $d-2$,
\item[(b)] for all $j = 3, \ldots, 5^{k}$, both $\mcG_1[V_i]$ and $\mcG_2[W_i]$ have exactly $r_{i,j}$
$j$-cycles, 
\item[(c)] for all  $j = 1, \ldots, 5^{k}$, both $\mcG_1[V_i]$ and $\mcG_2[W_i]$ have exactly $s_{i,j}$
$j$-paths with both endpoints of degree $d-1$, and
\item[(d)] there is an isomorphism 
$$f_0 : \mcG_1\Bigg[\bigcup_{i=1}^l NP(\mcG_1[V_i], 5^k)\Bigg] \ \to \ 
\mcG_2\Bigg[\bigcup_{i=1}^l NP(\mcG_2[W_i], 5^k)\Bigg]$$
such that, for all $p \in [l]$, if $v \in V_p$ then $f_0(v) \in W_p$.
\end{itemize}
\end{itemize}

\noindent
Clearly there is an integer $m$ such that if $\mcG \in \mbD$ and
$V_1, \ldots, V_l$ are the parts of the partition defined by $\xi$ in $\mcG$, then 
$\big| NP(\mcG[V_i], 5^k) \big| < \frac{m}{kl}$ for all $i$.
From the assumption that $c(\mbD) > 0$ and Theorem~\ref{unique decompositions}
it follows that the proportion of $\mcG \in \mbD \cap \mbP_n(l,d)$ such that
$\mcG$ has the $m$-extension property (with respect to the partition, defined
by $\xi$, on which its unique decomposition is based)
approaches 1 as $n \to \infty$.
Therefore, it suffices to prove that whenever both $\mcG_1, \mcG_2 \in \mbD \cap \mbP_n(l,d)$
have the $m$-extension property and $n$ is large enough, then $\mcG_1 \equiv_k \mcG_2$.

So suppose that both $\mcG_1, \mcG_2 \in \mbD \cap \mbP_n(l,d)$ have the $m$-extension property.
By~($*$), the parts of the partition of $V = [n]$ on which the unique decomposition of $\mcG_1$
is based can be ordered as $V_1, \ldots, V_l$ in such a way that (a)--(d) hold,
and the parts of the partition of $W = [n]$ on which the unique decomposition of $\mcG_2$
is based can be ordered as $W_1, \ldots, W_l$ in such a way that (a)--(d) hold.
In order to prove that $\mcG_1 \equiv_k \mcG_2$ it suffices to prove that
Duplicator has a winning strategy for the Ehrenfeucht-Fra\"{i}ss\'{e} game
in $k$ steps on $\mcG_1$ and $\mcG_2$. See for example \cite{EF} for 
definitions and results about Ehrenfeucht-Fra\"{i}ss\'{e} games.
We will use the following simplified notation instead of the one given
in Definition~\ref{definition of neighbourhood}. For $i = 1,2$ and $j = 1, \ldots, l$, if 
$v$ is a vertex of $\mcG_i[V_j]$, then $N(v, 5^k) = N_{\mcG_i[V_j]}(v, 5^k)$. 
To prove that Duplicator has a winning strategy for the 
Ehrenfeucht-Fra\"{i}ss\'{e} game on $\mcG_1$ and $\mcG_2$ it suffices to prove
the following statement.
\\

\noindent
{\em Claim.} Suppose that $i < k$, $v_1, \ldots, v_i \in V$, \
$w_1, \ldots, w_i \in W$, \
$v_{i+1} \in V$ (or $w_{i+1} \in W$)
and that 
\begin{align*}
f_i : \mcG_1\Bigg[\bigcup_{j=1}^l NP\big(\mcG_1[V_j], 5^{k-i}\big) \ \cup \ 
\bigcup_{j=1}^i &N\big(v_j, 5^{k-i}\big)\Bigg] \ \to \\ 
&\mcG_2\Bigg[\bigcup_{j=1}^l NP\big(\mcG_2[W_j], 5^{k-i}\big) \ \cup \
\bigcup_{j=1}^i N\big(w_j, 5^{k-i}\big)\Bigg]
\end{align*}
is an isomorphism such that $f_i$ extends $f_0$ from~(d), 
$f_i(v_j) = w_j$ for every $j = 1, \ldots, i$, and
whenever $v$ is in the domain of $f_i$ and $p \in [l]$, then $v \in V_p$ implies $f(v) \in W_p$.
Then there is $w_{i+1} \in W$ (or $v_{i+1} \in V$) such that all of the above
hold for `$i+1$' in place of `$i$'.
\\

\noindent
By symmetry it is enough to consider the case when $v_{i+1} \in V$ is given.
If 
\begin{equation}\label{neighbourhood of v-i+1 included in domain of f-i}
N\big(v_{i+1}, 5^{k-i-1}\big) \ \subseteq \ 
\bigcup_{j=1}^l NP\big(\mcG_1[V_j], 5^{k-i} - 1\big) \ \cup \ 
\bigcup_{j=1}^i N\big(v_j, 5^{k-i} - 1\big),
\end{equation}
then let $w_{i+1} = f(v_{i+1})$ and let $f_{i+1}$ be the restriction of $f_i$ to 
$$\bigcup_{j=1}^l NP\big(\mcG_1[V_j], 5^{k-i-1}\big) \ \cup \ \bigcup_{j=1}^i N\big(v_j, 5^{k-i-1}\big).$$
Now suppose that~(\ref{neighbourhood of v-i+1 included in domain of f-i}) does {\em not} hold,
so in particular $N\big(v_{i+1}, 5^{k-i-1}\big)$ is not included in 
$\bigcup_{j=1}^i N\big(v_j, 5^{k-i} - 1\big)$.
Let $p \in [l]$ be such that $v_{i+1} \in V_p$.
Then there is $u \in N\big(v_{i+1}, 5^{k-i-1}\big)$ such that, for every $j \leq i$, 
either $v_j \notin V_p$ or the distance,
in $\mcG_1[V_p]$, from $u$ to $v_j$ is at least $5^{k-i}$.
From this it follows that, for every $a \in N\big(v_{i+1}, 5^{k-i-1}\big)$
and every $j \leq i$, either $v_j \notin V_p$ or
the distance, in $\mcG_1[V_p]$, from $a$ to $v_j$ is at least
$$5^{k-i} - \mr{dist}_{\mcG_1[V_p]}(a,u) \ \geq \ 
5^{k-i} - 2\cdot 5^{k-i-1} \geq 3\cdot 5^{k-i-1}.$$
Consequently, for every $a \in N\big(v_{i+1}, 5^{k-i-1})$, every $j \leq i$
such that $v_j \in V_p$
and every $b \in N\big(v_j, 5^{k-i-1}\big)$, 
$$\dist_{\mcG_1[V_p]}(a,b) \ \geq \ 3 \cdot 5^{k-i-1} - 5^{k-i-1} \ = \ 2\cdot 5^{k-i-1} \ \geq \ 2,$$
because $i < k$. 
It follows that
\begin{align}\label{distance between v-i+1 and v-j}
&\text{if $j \leq i$ and $v_j \in V_p$, then the distance, in $\mcG_1[V_p]$, 
between } N\big(v_{i+1}, 5^{k-i-1}\big)\\
&\text{and $N\big(v_j, 5^{k-i-1}\big)$ is at least 2,} \nonumber
\end{align}
so, in particular, the two sets are disjoint.
In the same way, by considering any vertex in a small Poisson object
(Definition~\ref{definition of poisson object}) instead of $v_j$ for $j \leq i$,
it follows that 
\begin{equation}\label{distance between v-i+1 and poisson objects}
\text{the distance between $N\big(v_{i+1}, 5^{k-i-1}\big)$ and
$NP\big(\mcG_1[V_p], 5^{k-i-1}\big)$ is at least 2.}
\end{equation}
Now we are ready to use the assumption that $\mcG_1$ and $\mcG_2$ have the 
$m$-extension property. We start by defining an appropriate graph $\mcH$
as in Assumption~\ref{assumption about H}.
Let 
\begin{align*}
X_1 \ &= \ NP\big(\mcG_1[V_p], 5^{k-i-1}\big) \ \cup \
\bigcup \big\{ N\big(v_j, 5^{k-i-1}\big) : j \leq i \text{ and } v_j \in V_p \big\},\\
X_2 \ &= \ \bigcup \big\{ NP\big(\mcG_1[V_j], 5^{k-i-1}\big) : j \neq p \big\} \ \cup \ 
\bigcup \big\{ N\big(v_j, 5^{k-i-1}\big) : j \leq i \text{ and } v_j \notin V_p \big\},\\
Y \ &= \ N\big(v_{i+1}, 5^{k-i-1}\big).
\end{align*}
Then let $\mcH = \mcG_1[X_1 \cup X_2 \cup Y]$.
By~(\ref{distance between v-i+1 and v-j})
and~(\ref{distance between v-i+1 and poisson objects}),
$X_1$, $X_2$ and $Y$ are mutually disjoint.
By assumption, $f_i \uhrc X_1 \cup X_2$ is a strong embedding of $\mcH[X_1 \cup X_2]$
into $\mcG_2$ such that $f(X_1) \subseteq W_p$ and $f(X_2) \subseteq W \setminus W_p$.
From~(\ref{distance between v-i+1 and poisson objects}) and
the assumption that $\mcG_1 \in \mbD \subseteq \mbP^k(l,d)$ it follows that
the subgraph of $\mcG_1[V_p]$ induced by $Y = N(v_{i+1}, 5^{k-i-1})$ is characterised as follows:
\begin{itemize}
\item[] {\em Case $d = 2$:} A path with at least $5^{k-i-1}$ vertices in which the distance
from $v_{i+1}$ to at least one of the endpoints is $5^{k-i-1}$.

\item[] {\em Case $d \geq 3$:} A tree in which every path from $v_{i+1}$ to a leaf has length 
$5^{k-i-1}$ and either every non-leaf has degree $d$, 
or exactly one non-leaf has degree $d-1$ and all other non-leaves have degree $d$.
\end{itemize}
Since $\mcG_2 \in \mbD \subseteq \mbP^k$ it follows that there are
between $\sqrt{(d-\varepsilon)n}$ and $\sqrt{(d+\varepsilon)n}$ vertices in $\mcG_2[V_p]$
with degree $d-1$ (in $\mcG_2[V_p]$). Moreover, because of the properties of $f_i$ and $f_0$,
for any two distinct vertices $w, w' \in W_p \setminus f_i(X_1 \cup X_2)$ with degree $d-1$
in $\mcG_2[W_p]$ the distance between them is at least $5^{k+1} + 1$.
Hence, for every $w \in W_p \setminus f_i(X_1 \cup X_2)$, the subgraph of $\mcG_2[W_p]$ induced
by $N(w, 5^{k-i-1})$ is characterised in the same way as in the case of
$Y = N(v_{i+1}, 5^{k-i-1})$, as explained above.
Therefore, there are (assuming $|W|$ is large enough) at least $n^{1/4}$ different
induced subgraphs of $\mcG_2[W_p]$ which are isomorphic to $\mcH[Y]$.
By the choice of $m$ we have $|X_1 \cup X_2 \cup Y| \leq m$ and since
$\mcG_2$ has the $m$-extension property there is
a strong embedding $f_{i+1}$ of $\mcH$ into $\mcG_2$ which extends $f_i$
and such that $f(Y) \subseteq W_p$.
If we let $w_{i+1} = f_{i+1}(v_{i+1})$ then $f_{i+1}(v_j) = w_j$ for every $j = 1, \ldots, i+1$,
and whenever $v$ is in the domain of $f_{i+1}$ and $j \in [l]$, then
$v \in V_j$ implies $f_{i+1}(v) \in W_j$.
\hfill $\square$
\\

\noindent
Let $(\mbD_i : i \in \mbbN)$ be an enumeration of all $\approx_k^+$-classes
and let, by Corollary~\ref{convergence of strong poisson classes}, $0 \leq d_i < 1$ be the constant
$$d_i = \lim_{n\to\infty} \frac{|\mbD_i \cap \mbP_n(l,d)|}{|\mbP_n(l,d)|}.$$
By the same corollary we have 
\begin{equation}\label{the sum of all d-i is 1}
\sum_{i=0}^{\infty} d_i \ = \  1.
\end{equation}
Also note that since, for every $n$, the equivalence relation $\approx_k^+$ partitions the (finite) set 
$\mbP_n(l,d)$, it follows that
\begin{equation}\label{all classes sum to one}
\text{for every $n$ and every $\mbA \subseteq \mbP_n(l,d)$, } \ \ 
\frac{|\mbA \cap \mbP_n(l,d)|}{|\mbP_n(l,d)|} \ =  \ \sum_{i=0}^{\infty} \frac{|\mbA \cap \mbD_i \cap \mbP_n(l,d)|}{|\mbP_n(l,d)|}.
\end{equation}

\begin{cor}\label{convergence of any subset of classes}
Let $\mbE$ be any $\equiv_k$-class.\\
(i) For every $i \in \mbbN$,
$$\lim_{n\to\infty} \ \frac{|\mbE \cap \mbD_i \cap \mbP_n(l,d)|}{|\mbP_n(l,d)|} \ \ 
\text{ equals either $0$ or $d_i$.}$$
(ii) Let $I$ be the set of $i \in \mbbN$ such that the above limit equals $d_i$.
Then 
$$\lim_{n\to\infty} \ \frac{|\mbE \cap \mbP_n(l,d)|}{|\mbP_n(l,d)|} \ = \ 
\sum_{i \in I} d_i.$$
\end{cor}

\noindent
{\bf Proof.}
Part (i) is a consequence, using the notation `$d_i$', of 
Lemma~\ref{convergence of k-equivalence class cut with strong poisson class},
so we turn to the proof of~(ii).
We use the abbreviation $\mbP_n = \mbP_n(l,d)$.
Let $\varepsilon > 0$. We show that for $I$ as defined above and large enough $n$,
$$\Bigg|\frac{|\mbE \cap \mbP_n|}{|\mbP_n|} \ - \ \sum_{i \in I} d_i\Bigg| 
\ \leq \ 5\varepsilon.$$
By (\ref{the sum of all d-i is 1}) we can choose $m$ large enough that 
\begin{equation}\label{the choice of m}
\sum_{i = m+1}^{\infty} d_i \ \leq \varepsilon \quad \text{ and } \quad
1 - \varepsilon \ \leq \ \sum_{i=0}^m d_i \ \leq 1.
\end{equation} 
Then, using the definition of $d_i$, $i \in \mbbN$, and part~(i),
choose $n_0$ so that
\begin{equation}\label{choice of n-0, part one}
\text{for all $i \leq m$ such that $i \notin I$ and all $n > n_0$}, \ 
\frac{|\mbE \cap \mbD_i \cap \mbP_n|}{|\mbP_n|} \ \leq \ \frac{\varepsilon}{m+1},
\end{equation}
\begin{equation}\label{choice of n-0, part two}
\text{for all $i \leq m$ such that $i \in I$ and all $n > n_0$}, \ 
\Bigg| \frac{|\mbE \cap \mbD_i \cap \mbP_n|}{|\mbP_n|} \ - \ d_i \Bigg| \ \leq \ \frac{\varepsilon}{m+1},
\end{equation}
and
\begin{equation}\label{choice of n-0, part three}
\text{for all $i \leq m$ and all $n > n_0$}, \ 
\Bigg| \frac{|\mbD_i \cap \mbP_n|}{|\mbP_n|} \ - \ d_i \Bigg| \ \leq \ \frac{\varepsilon}{m+1},
\end{equation}
For all $n > n_0$ we get, by the use 
of~(\ref{all classes sum to one})--(\ref{choice of n-0, part three}),
\begin{align*}
&\Bigg| \frac{|\mbE \cap \mbP_n|}{|\mbP_n|} \ - \ \sum_{i \in I} d_i \Bigg| \ = \ 
\Bigg| \sum_{i=0}^{\infty} \frac{|\mbE \cap \mbD_i \cap \mbP_n|}{|\mbP_n|} \ - \ \sum_{i \in I} d_i \Bigg| \\
\leq \ 
&\sum_{i\leq m, i\notin I} \frac{|\mbE \cap \mbD_i \cap \mbP_n|}{|\mbP_n|} \ + \ 
\sum_{i\leq m, i\in I} \Bigg| \frac{|\mbE \cap \mbD_i \cap \mbP_n|}{|\mbP_n|} \ - \ d_i \Bigg| \\ 
&+ \ 
\Bigg| \sum_{i=m+1}^{\infty} \frac{|\mbE \cap \mbD_i \cap \mbP_n|}{|\mbP_n|} \Bigg| \ + \ 
\Bigg| \sum_{i=m+1}^{\infty} d_i \Bigg| \\
\leq \ &3\varepsilon \ + \ \Bigg| \sum_{i=m+1}^{\infty} \frac{|\mbD_i \cap \mbP_n|}{|\mbP_n|} \Bigg| 
\ = \ 
3\varepsilon \ + \ \Bigg| \sum_{i=0}^{\infty} \frac{|\mbD_i \cap \mbP_n|}{|\mbP_n|} \ - \ 
\sum_{i\leq m} \frac{|\mbD_i \cap \mbP_n|}{|\mbP_n|} \Bigg| \\
= \ &3\varepsilon \ + \ \Bigg| 1 \ - \ \sum_{i\leq m} \frac{|\mbD_i \cap \mbP_n|}{|\mbP_n|} \Bigg|
\ \leq \ 3\varepsilon \ + \ 
\Bigg| \varepsilon \ + \ \sum_{i \leq m} d_i \ - \ \sum_{i\leq m} \frac{|\mbD_i \cap \mbP_n|}{|\mbP_n|} \Bigg|
\\
\leq \ &4\varepsilon \ + \ 
\sum_{i \leq m} \Bigg| d_i \ - \ \frac{|\mbD_i \cap \mbP_n|}{|\mbP_n|} \Bigg| \ \leq \ 
5\varepsilon. \hspace{40mm} \square
\end{align*}

\noindent
Corollary~\ref{convergence of any subset of classes} concludes the proof of
the limit law of $\mbP(l,d)$.
The second claim of part~(ii) of Theorem~\ref{limit theorem},
that $\mbP(l,d)$ does not have a zero-one law, follows from the fact given
by Lemma~\ref{convergence of poisson classes} that (for any $k \geq 1$) one can choose infinitely many
$\approx_k$-classes $\mbC_i$ and sentences $\varphi_i$, $i \in \mbbN$, 
such that $|\mbC_i \cap \mbP_n(l,d)| \big/ |\mbP_n(l,d)|$
converges to a positive number, and $\varphi_i$ is true in every $\mcG \in \mbC_i$ and
false in every $\mcG \in \mbC_j$ if $j \neq i$.
For example, let $\mbC_i$ be an $\approx_k$-class such that for every $\mcG \in \mbC_i$
and every part $V_j$ of the partition defined by $\xi(x,y)$, $\mcG[V_j]$ has exactly
$i$ vertices with degree $d-2$, 
and let $\varphi_i$ express this.

\subsection{Proof of part (i) of Theorem~\ref{limit theorem}}
\label{proof of part one of limit law for almost l-partite graphs}

\noindent
We now sketch the proof of part~(i) of Theorem~\ref{limit theorem}.
By Lemma~2.11 in \cite{Kop12},
for every $\varepsilon > 0$,
the proportion of $\mcG \in \mbP_n(1,1)$ such that
there are between $n^{1/2 - \varepsilon}$ and $n^{1/2 + \varepsilon}$ vertices
with degree 0, approaches 1 as $n \to \infty$.
This information about $\mbP(1,1)$ is sufficient for proving a zero-one law for $\mbP(l,1)$.

Fix $0 < \varepsilon < 1/4$.
By similar reasoning as when proving Lemma~\ref{P-t is almost everything}
it follows that, for every $k \in \mbbN$, the proportion $\mcG \in \mbP_n(l,1)$ with the following
properties approaches 1 as $n \to \infty$:
\begin{itemize}
\item[(a)] $\mcG$ has a unique decomposition which is based on a $\mu n$-rich partition
$V_1, \ldots, V_l$ of $[n]$ which is defined by $\xi(x,y)$ (from Definition~\ref{definition of xi}), 
\item[(b)] $\mcG$ has the $2k$-extension property with respect to the partition defined by $\xi(x,y)$, and
\item[(c)] for every $i \in [l]$, $\mcG[V_i]$ has between $n^{1/2 - \varepsilon}$
and $n^{1/2 + \varepsilon}$ vertices degree 0.
\end{itemize}
Hence, it suffices to prove, for an arbitrary integer $k > 0$, 
that if $\mcG_1, \mcG_2 \in \mbP_n(l,1)$ satisfy~(a),~(b) and~(c),
then $\mcG_1 \equiv_k \mcG_2$.
This is done in a similar way as we proved, in the proof of
Lemma~\ref{convergence of k-equivalence class cut with strong poisson class},
that if $\mcG_1, \mcG_2 \in \mbD$ have the $m$-extension property, for suitably chosen $m$,
then $\mcG_1 \equiv_k \mcG_2$.

\section{Forbidden subgraphs}\label{forbidden subgraphs}

\noindent
Recall that for integers $1 \leq s_1 \leq \ldots \leq s_l$, $\mcK_{1, s_1, \ldots, s_l}$ denotes the 
complete $l$-partite graph with parts (or colour classes) of sizes $1, s_1, \ldots, s_l$.
For any graph $\mcH$, $\mb{Forb}_n(\mcH)$ denotes the set of graphs with vertices
$1, \ldots, n$ in which there is no subgraph isomorphic to $\mcH$ and 
$\mb{Forb}(\mcH) = \bigcup_{n\in \mbbN^+} \mb{Forb}_n(\mcH)$.
In this section we use Theorems~\ref{HPS-results} and~\ref{limit theorem} together with
some new technical results to prove:
\medskip

\noindent
{\bf Theorem~\ref{limit law for forbidden l+1-partite graphs}.} 
{\em
Suppose that $l \geq 2$, $1 \leq s_1 \leq s_2 \leq \ldots \leq s_l$ are integers.\\
(i) For every sentence $\varphi$ in the language of graphs, the 
proportion of $\mcG \in \mb{Forb}_n(\mcK_{1,s_1, \ldots, s_l})$ in which $\varphi$ is true
converges as $n \to \infty$. \\
(ii) If $s_1 \leq 2$ then this proportion converges to 0 or 1 for every sentence $\varphi$.\\
(iii) If $s_1 > 2$ then there are infinitely many mutually contradictory sentences 
$\varphi_i$, $i \in \mbbN$, in the language of graphs such that the proportion of
$\mcG \in \mb{Forb}_n(\mcK_{1,s_1, \ldots, s_l})$ in which $\varphi_i$ is true approaches
some $\alpha_i$ such that $0 < \alpha_i < 1$.
}
\medskip

\noindent
Part (i) of Theorem~\ref{limit law for  forbidden l+1-partite graphs} 
is an immediate consequence of Lemmas~\ref{existence of a cycle}
and~\ref{the convergence lemma for H} below.
Part~(ii) and~(iii) of Theorem~\ref{limit law for  forbidden l+1-partite graphs} 
require some more argumentation, which
is given after the proof of Lemma~\ref{the convergence lemma for H}.

\begin{lem}\label{existence of a cycle}
Suppose that $l \geq 1$ and $1 \leq s_1 \leq \ldots \leq s_l$ are integers.
If $V_1, \ldots, V_l$ is a partition of the vertex set of $\mcK_{1, s_1, \ldots, s_l}$
such that, for every $i = 1, \ldots, l$ and every $v \in V_i$,
$v$ has at most $s_1 - 1$ neighbours in $V_i$,
then, for some $i$, $V_i$ contains a 3-cycle.
\end{lem}

\noindent
{\bf Proof.}
Let $C_0, \ldots, C_l$ be the ``colour classes'' of $\mcK_{1, s_1, \ldots, s_l}$,
in other words, $\bigcup_{i = 0}^l C_i$ is the vertex set of $\mcK_{1, s_1, \ldots, s_l}$
and vertices $v$ and $w$ are adjacent to each other if and only if
there are $i \neq j$ such that $v \in V_i$ and $w \in V_j$.
Moreover, assume that $|C_0| = 1$ and $|C_i| = s_i$ for $i = 1, \ldots, l$.
We use induction on $l$.
It is clear that $\mcK_{1, s_1}$ cannot be ``partitioned'' into one class such that
every vertex has at most $s_1 - 1$ neighbours in its own part.
This takes care of the base case $l = 1$.

Now assume that $l > 1$ and that $V_1, \ldots, V_l$ is a 
partition of $\bigcup_{i = 0}^l C_i$ such that for every $i = 1, \ldots, l$
and every $v \in V_i$, $v$ has at most $s_1 - 1$ neighbours in $V_i$.
As $C_0$ is a singleton, we have $C_0 \subseteq V_i$ for some $i$.
By reordering $V_1, \ldots, V_l$ if necessary, we may assume that $i = 1$.
If there are $i > j > 0$ such that $V_1 \cap C_i \neq \es$ and $V_1 \cap C_j \neq \es$,
then, as $C_0 \subseteq V_1$, it follows that $V_1$ contains a 3-cycle and we are done.
So now suppose that there is at most one $i > 0$ such that $V_1 \cap C_i \neq \es$.
First suppose that there exists exactly one $k > 0$ such that $V_1 \cap C_k \neq \es$.
Since $C_0 \subseteq V_1$ and every vertex in $V_1$ has at most $s_1 - 1$ neighbours in $V_1$
it follows that $|V_1 \cap C_k| \leq s_1 - 1 \leq s_k - 1$
which implies that $C_k \setminus V_1 \neq \es$.
Choose any vertex $v \in C_k \setminus V_1$.
Consider the subgraph $\mcK'$ of $\mcK_{1, s_1, \ldots, s_l}$ which is induced
by 
$$C' = \{v\} \cup \bigcup_{1 \leq i \leq l, i \neq k} C_k,$$
so in model theoretic terms, $\mcK'$ is the substructure of $\mcK_{1, s_1, \ldots, s_l}$ with 
universe $C'$.
Note that $\mcK'$ is a complete $l$-partite graph with $l$ colour classes
$\{v\}, C_1, \ldots, C_{k-1}, C_{k+1}, \ldots, C_l$, one of which is a singleton.
Also, $V_2 \cap C', \ldots, V_l \cap C'$ is a partition of $C'$ such that for
every $i = 2, \ldots, l$ and every $w \in V_i$, $w$ has at most $s_1 - 1$ neighbours in $V_i$.
Since $s_1 \leq s_2 \leq \ldots \leq s_l$ are the cardinalities of 
$C_1, \ldots, C_l$, respectively, it follows that, for every $i = 1, \ldots, l$, $|C_i| \geq s_1$.
Therefore the induction hypothesis implies that for some $i \geq 2$,
$V_i \cap C'$ contains a 3-cycle, and the same 3-cycle is contained in $V_i$.

Now suppose that $C_i \cap V_1 = \es$ for every $i > 0$, so $V_1 = C_0$ (by the assumption
that $C_0 \subseteq V_1$).
Then we can choose any $v \in C_1$ and consider the subgraph $\mcK'$ of 
$\mcK_{1, s_1, \ldots, s_l}$ which is induced by the set
$C' = \{v\} \cup \bigcup_{i = 2}^l C_i$. Note that $V_2 \cap C', \ldots, V_l \cap C'$
is a partition of $C'$ such that for all $i = 2, \ldots, l$
and every $v \in V_i \cap C'$, $v$ has at most $s_1 - 1 \leq s_2 - 1$ neighbours
in $V_i$. Therefore the induction hypothesis implies that for some $i \geq 2$,
$V_i \cap C'$ contains a 3-cycle.
\hfill $\square$

\noindent
Recall the definition of a colour-critical vertex and criticality of such a vertex,
defined before the statement of Theorem~\ref{HPS-results}.

\begin{lem}\label{the convergence lemma for H}
Suppose that $l \geq 2$, $d \geq 0$ and that $\mcH$ is a graph with the following properties:
\begin{itemize}
\item[\textbullet] \ $\mcH$ has chromatic number $l+1$.
\item[\textbullet] \ $\mcH$ has a colour-critical vertex with criticality $d+1$ and
no colour-critical vertex has criticality smaller than $d+1$.
\item[\textbullet] \ If $V_1, \ldots, V_l$ is a partition of the vertex set of $\mcH$
such that for, every $i = 1, \ldots, l$, and every $v \in V_i$, $v$ has at most
$d$ neighbours in $V_i$, then, for some $j$, $V_j$ contains a 3-cycle.
\end{itemize}
Then, for every sentence $\varphi$ in the language of graphs, the proportion
of graphs $\mcG \in \mb{Forb}_n(\mcH)$ such that $\mcG \models \varphi$ converges as $n \to \infty$.
\end{lem}

\noindent
{\bf Proof.}
Let $\mcH$ be a graph with the listed properties.
Let $\mbF_n(\mcH) = \mb{Forb}_n(\mcH)$, $\mbP_n = \mbP_n(l,d)$ and let
$\mbX_n$ be the set of $\mcG \in \mbP_n$ such that
$\mcG$ has a unique decomposition based on a partition $V_1, \ldots, V_l$ definable by $\xi(x,y)$, as in 
Theorem~\ref{unique decompositions}, and for every $i = 1, \ldots, l$,
$V_i$ contains {\em no} 3-cycle.
Observe that it follows from the third property of $\mcH$ (listed in the lemma)
that $\mbX_n \subseteq \mb{Forb}_n(\mcH)$ for all $n$.

Until further notice, assume that $d \geq 2$ 
(as in Section~\ref{proof of part two of limit law for almost l-partite graphs})
Choose any $k \geq 1$.
Recall Definition~\ref{definition of poisson equivalence} of the equivalence relation `$\approx_k$'.
Let $\mbC \subseteq \mbP(l,d)$ be the union of all $\approx_k$-classes included in
$\bigcup_{n\in \mbbN^+}\mbP_n^k(l,d)$ 
(see Definition~\ref{definition of P-t}) and such that
if $\mcG \in \mbC$, then {\em no} part of the partition of the vertex set of
$\mcG$ defined by $\xi(x,y)$ contains a 3-cycle.
By Lemma~\ref{convergence of poisson classes}, there is $c > 0$
such that 
\begin{equation}\label{convergence of graphs without 3-cycles}
\lim_{n\to\infty} \ \frac{\big|\mbC \cap \mbP_n\big|}{\big|\mbP_n\big|} \ = \ c.
\end{equation}
By Theorem~\ref{unique decompositions}, the proportion of $\mcG \in \mbP_n$ which have
a unique decomposition and the partition on which it is based is defined by $\xi(x,y)$,
approaches 1 as $n \to \infty$. From this,~(\ref{convergence of graphs without 3-cycles})
and since $c > 0$ it follows that 
\begin{equation}\label{convergence of graphs without 3-cycles and with unique decomposition}
\lim_{n\to\infty} \ \frac{\big|\mbX_n|}{\big|\mbP_n\big|} \ = \ c.
\end{equation}
Let $\psi_{\mcH}$ be a sentence which expresses that ``there is no subgraph 
isomorphic to $\mcH$'', so for every graph $\mcG$ with vertices $1, \ldots, n$, 
$\mcG \in \mb{Forb}_n(\mcH)$ if and only if $\mcG \models \psi_{\mcH}$.
Then $\mbF_n(\mcH) \cap \mbP_n = \big\{\mcG \in \mbP_n : \mcG \models \psi_{\mcH}\big\}$
and from Theorem~\ref{limit theorem}
it follows that for some $0 \leq b \leq 1$,
\begin{equation}\label{convergence of H-free graphs in P-n}
\lim_{n\to\infty} \ \frac{\big|\mbF_n(\mcH) \cap \mbP_n|}{\big|\mbP_n\big|} \ = \ 
\frac{\big|\{\mcG \in \mbP_n : \mcG \models \psi_{\mcH}\}\big|}{\big|\mbP_n\big|} \ = \ b.
\end{equation}
Since $c > 0$ and $\mbX_n \subseteq \mbF_n(\mcH) \cap \mbP_n$,
it follows from~(\ref{convergence of graphs without 3-cycles and with unique decomposition})
and~(\ref{convergence of H-free graphs in P-n})
that $b \geq c > 0$.
We have arrived at this conclusion under the assumption that $d \geq 2$.
If $d = 0$ or $d = 1$ and the vertex set of $\mcG \in \mbP_n(l,d)$ is partitioned
into $l$ parts $V_1, \ldots, V_l$ such that no vertex has more than $d$ neighbours in its own part,
then $V_i$ clearly does not contain a 3-cycle for any $i$.
So if $d = 0$ or $d = 1$, then $|\mbX_n| \big/ |\mbP_n(l,d)|$ converges to 1 as $n \to \infty$,
by Theorem~\ref{unique decompositions}.
Hence we get~(\ref{convergence of H-free graphs in P-n}) for some $b > 0$ also in the case $d \in \{0,1\}$.
The rest of the proof does not depend on whether $d \leq 1$ or $d \geq 2$.

Let $\varphi$ be any sentence in the language of graphs.
Then, for large enough $n$,
\begin{align*}
&\frac{\big|\{ \mcG \in \mbF_n(\mcH) : \mcG \models \varphi\}\big|}{\big|\mbF_n(\mcH)\big|} \\ 
= \ 
&\frac{\big|\{  \mcG \in \mbF_n(\mcH) \cap \mbP_n : \mcG \models \varphi \}\big|}{\big|\mbF_n(\mcH)\big|}
\ + \ 
\frac{\big|\{ \mcG \in \mbF_n(\mcH) \setminus \mbP_n : \mcG \models \varphi \}\big|}{\big|\mbF_n(\mcH)\big|} \\
= \
&\frac{\big|\{ \mcG \in \mbP_n : \mcG \models \psi_{\mcH} \wedge \varphi \}\big|}{\big|\mbF_n(\mcH)\big|} 
\ + \ 
\frac{\big|\{ \mcG \in \mbF_n(\mcH) \setminus \mbP_n : \mcG \models \varphi \}\big|}{\big|\mbF_n(\mcH)\big|} \\
= \ 
&\frac{\big|\{ \mcG \in \mbP_n : \mcG \models \psi_{\mcH} \wedge \varphi \}\big|}
{\big|\mbP_n\big|}
\ \cdot \ 
\frac{\big|\mbP_n\big|}{\big|\mbF_n(\mcH) \cap \mbP_n\big|}
\ \cdot \ 
\frac{\big|\mbF_n(\mcH) \cap \mbP_n\big|}{\big|\mbF_n(\mcH)\big|} \\
&+ \ 
\frac{\big|\{ \mcG \in \mbF_n(\mcH) \setminus \mbP_n : \mcG \models \varphi \}\big|}{\big|\mbF_n(\mcH)\big|} \\
\to \ 
&\Bigg( \lim_{n\to\infty} \ \frac{\big|\{ \mcG \in \mbP_n : \mcG \models \psi_{\mcH} \wedge \varphi \}\big|}
{\big|\mbP_n\big|} \Bigg) \ \cdot \ \frac{1}{b} \quad \text{ as } n \to \infty,
\end{align*}
because of~(\ref{convergence of H-free graphs in P-n}), where $b > 0$ as explained above,
and Theorems~\ref{HPS-results} and~\ref{limit theorem}.
\hfill $\square$
\\

\noindent
Part~(i) of Theorem~\ref{limit law for  forbidden l+1-partite graphs}
follows directly from Lemmas~\ref{existence of a cycle}
and~\ref{the convergence lemma for H}.
Now we consider part~(ii) of Theorem~\ref{limit law for  forbidden l+1-partite graphs}.
Suppose that $l \geq 2$ and $1 \leq s_1 \leq \ldots \leq s_l$ are integers and
that $s_1 \leq 2$.
Then $\mcH = \mcK_{1, s_1, \ldots, s_l}$ has the three properties listed
in Lemma~\ref{the convergence lemma for H}. 
The proof of Lemma~\ref{the convergence lemma for H}
shows (also in the case $s_1 \leq 2$)
that the proportion of $\mcG \in \mbP_n(l, s_1 - 1)$ which are 
$\mcK_{1, s_1, \ldots, s_l}$-free converges to a positive number.
From Theorem~\ref{HPS-results} and 
Theorem~\ref{limit theorem}~(i)
it follows that, for every sentence $\varphi$, the proportion of $\mcG \in \mb{Forb}_n(\mcH)$
in which $\varphi$ is true converges to either 0 or 1.

Now we prove part~(iii) of Theorem~\ref{limit law for  forbidden l+1-partite graphs}.
Suppose that $s_1 \geq 3$.
Choose any integer $k \geq 1$ and note that $5^k > 4$. 
Then pick any integers $p > q \geq 0$.
We argue similarly as in the proof of Lemma~\ref{the convergence lemma for H}.
Let $\mbC, \mbD \subseteq \mbP(l, s_1 - 1)$ be $\approx_k$-equivalence classes 
(see Definition~\ref{definition of poisson equivalence}) such that 
$\mbC, \mbD \subseteq \bigcup_{n \in \mbbN^+} \mbP_n^k(l, s_1 - 1)$ and the following hold:
\begin{itemize}
\item[(a)] If $\mcG \in \mbC$ and $V_1, \ldots, V_l$ is the partition of the vertex
set of $\mcG$ which is defined by $\xi(x,y)$, then for every $i \in [l]$,
$\mcG[V_i]$ has no 3-cycle and exactly $p$ vertices with degree $s_1 - 3$.
\item[(b)] If $\mcG \in \mbD$ and $V_1, \ldots, V_l$ is the partition of the vertex
set of $\mcG$ which is defined by $\xi(x,y)$, then for every $i \in [l]$,
$\mcG[V_i]$ has no 3-cycle and exactly $q$ vertices with degree $s_1 - 3$.
\end{itemize}
Note that $\mbC$ and $\mbD$ are distinct
$\approx_k$-equivalence classes, so $\mbC \cap \mbD = \es$.
By Lemma~\ref{convergence of poisson classes},
$|\mbC \cap \mbP_n(l, s_1 - 1)| \big/ |\mbP_n(l, s_1 - 1)|$ converges to a positive number
as $n \to \infty$,
and the same holds for $\mbD$ in place of $\mbC$.
By Theorem~\ref{unique decompositions},
the proportion of graphs $\mcG \in \mbC \cap \mbP_n(l, s_1 - 1)$ 
(respectively $\mcG \in \mbD \cap \mbP_n(l, s_1, - 1)$)
such that $\mcG$ has a unique decomposition and this decomposition 
is based on a partition defined by $\xi(x,y)$, approaches 1 as $n \to \infty$.
It follows, by using Lemma~\ref{existence of a cycle}, that the
proportion of $\mcG \in \mbC \cap \mbP_n(l, s_1 - 1)$ 
(respectively $\mcG \in \mbD \cap \mbP_n(l, s_1 - 1)$)
that are $\mcK_{1, s_1, \ldots, s_l}$-free approaches 1 as $n \to \infty$.
In the proof of Lemma~\ref{the convergence lemma for H},
which is applicable to $\mcH = \mcK_{1, s_1, \ldots, s_l}$ (by Lemma~\ref{existence of a cycle}),
it was shown, see~(\ref{convergence of H-free graphs in P-n}), that
$$\frac{\big|\mbP_n(l, s_1 - 1) \cap \mb{Forb}_n(\mcK_{1, s_1, \ldots, s_l})\big|}
{\big|\mbP_n(l, s_1 - 1)\big|}
\quad \text{ converges to a positive number.}$$
These conclusions together with Theorem~\ref{HPS-results}
imply that both the quotients
$$\frac{\big|\mbC \cap \mb{Forb}_n(\mcK_{1, s_1, \ldots, s_l})\big|}
{\big|\mb{Forb}_n(\mcK_{1, s_1, \ldots, s_l})\big|} 
\quad \text{and} \quad
\frac{\big|\mbD \cap \mb{Forb}_n(\mcK_{1, s_1, \ldots, s_l})\big|}
{\big|\mb{Forb}_n(\mcK_{1, s_1, \ldots, s_l})\big|}$$
converge to positive numbers.
Since $\mbC \cap \mbD = \es$ it follows that none of these numbers can be~1.
As $p > q \geq 0$ where arbitrary and the property ``the induced subgraph 
on every part of the partition defined by $\xi(x,y)$ has no 3-cycle and exactly
$p$ vertices with degree $s_1 - 3$'' can be expressed with the (first-order) language of graphs
this completes the proof of part~(iii) of Theorem~\ref{limit law for  forbidden l+1-partite graphs},
and hence the proof of that theorem is finished.

\begin{rem}{\rm
Suppose that $l \geq 2$ and $1 \leq s_1 \leq \ldots \leq s_l$ are integers.
One can prove, by a combinatorial argument, that `$s_1 \leq 2$ or $s_2 \geq 2(s_1 - 1)$'
is a necessary and sufficient condition for $\mcK_{1, s_1, \ldots, s_l}$ having the property:
there is a partition of the vertex set such that every vertex has at most $s_1 - 1$ 
neighbours in its own part.
Consequently, $\mbP(l, s_1 - 1) \subseteq \mb{Forb}_n(\mcK_{1, s_1, \ldots, s_l})$
if and only if $s_1 \leq 2$ or $s_2 \geq 2(s_1 - 1)$.
It also follows that 
$|\mbP_n(l, s_1 - 1) \cap \mb{Forb}_n(\mcK_{1, s_1, \ldots, s_l})| \big/ 
|\mbP_n(l, s_1 - 1)|$ converges to 1, as $n \to \infty$, if and only if 
$s_1 \leq 2$ or $s_2 \geq 2(s_1 - 1)$; otherwise this ratio converges to a positive number less than 1.
}\end{rem}

\end{document}